\documentclass[10pt,a4paper,oneside,leqno]{article}
\usepackage{amsfonts,amssymb}
\usepackage{eufrak}
\usepackage{amscd}
\usepackage{color}

\newtheorem{defn}{Definition}
\newtheorem{theoremn}{Theorem}
\newtheorem{observen}{Observation}
\newtheorem{example}{Example}

\vspace{0.5cm}

 4
\font\ebf=cmbx8
\font\erm=cmr8

\usepackage{graphicx, graphics}

\begin{document}

% % % New Commands

\newcommand{\fnomial}[2]{ {{#1} \choose {#2}}_{\!\!F} }
\newcommand{\fnomialF}[3]{ {{#1} \choose {#2}}_{\!\!#3} }

\thispagestyle{empty}

\noindent {\bf Cobweb posets - Recent Results}

\vspace{0.3cm}

\noindent A. Krzysztof Kwa\'sniewski (*) ,M.~Dziemia{\'n}czuk (**)

\vspace{0.15cm}
\noindent {\erm (*) the Dissident - relegated by Bia\l ystok University authorities  }\\
\noindent {\erm from the Institute of Informatics to Faculty of Physics}\\
\noindent {\erm ul. Lipowa 41,  15 424  Bia\l ystok, Poland}\\
\noindent {\erm e-mail: kwandr@gmail.com}\\
\noindent {\erm (**) Former  Student in the Institute of Computer Science, Bia\l ystok University; then moved to:} \\
\noindent {\erm Institute of Informatics, University of Gda\'nsk, Poland; }\\
\noindent {\erm e-mail: mdziemianczuk@gmail.com} 

\vspace{0.4cm}

\noindent {\ebf SUMMARY}

\vspace{0.1cm}

\noindent {\small  Cobweb posets uniquely represented by directed acyclic graphs are such a generalization of the Fibonacci tree that allows joint combinatorial interpretation for all of them under admissibility condition. This interpretation was derived in the source papers ([6,7] and references therein to the first author).[7,6,8] include  natural enquires to be reported on here. The purpose of this presentation is to report on the progress in solving  computational problems which are quite easily formulated for the new class of directed acyclic graphs
interpreted as Hasse diagrams. The problems posed there and  not yet all solved completely are of crucial 
importance  for the vast class of  new partially ordered sets with joint combinatorial interpretation.
These so called cobweb posets - are relatives of Fibonacci tree and are labeled by  specific
number sequences - natural numbers sequence and Fibonacci sequence included.  One presents here also a  
join combinatorial interpretation of those posets` $F$-nomial coefficients which are computed
with the so called cobweb admissible sequences. Cobweb posets and their natural subposets are 
graded posets. They are vertex partitioned into such antichains $\Phi_n$ 
(where $n$ is a nonnegative integer) that for each $\Phi_n$, all of the elements
covering $x$ are in $\Phi_{n+1}$ and all the elements covered by $x$ are in $\Phi_n$. We shall call
the $\Phi_n$ the  $n-th$- level. The cobweb posets might  be identified with a chain of di-bicliques i.e. 
by definition - a chain of complete bipartite one direction digraphs [6].
Any chain of relations is therefore obtainable from the cobweb poset chain of complete relations
via deleting arcs in di-bicliques of the complete relations chain.
In particular  we response to one of those problems  [1].  This is a tiling problem. Our information on tiling problem refers on proofs of tiling's existence for some cobweb-admissible sequences as in [1]. There the second author shows that not all cobwebs admit tiling as defined below and 
provides examples of  cobwebs admitting tiling.}

\vspace{0.3cm}

\noindent Key Words: acyclic digraphs, tilings, special number sequences, binomial-like coefficients.

\vspace{0.1cm}

\noindent AMS Classification Numbers: 06A07 ,05C70,05C75, 11B39.

\vspace{0.1cm}

\noindent  Presented at Gian-Carlo Polish Seminar:

\noindent \emph{http://ii.uwb.edu.pl/akk/sem/sem\_rota.htm}

\vspace{0.1cm}

\noindent published:  Adv. Stud. Contemp. Math.  vol. \textbf{16} (2)  April (\textbf{2008}):197-218. 

\vspace{0.4cm}

\section{Introduction [6]}
\noindent A directed acyclic graph, also called DAG, is a directed graph with no directed cycles. 
DAGs  considered  as a generalization of trees  have a lot of applications in computer science,
bioinformatics,  physics and many natural activities of humanity and nature. 
Here we introduce specific DAGs as generalization of trees being
inspired by algorithm of the Fibonacci tree growth. For any given natural numbers valued sequence
the graded (layered) cobweb posets` DAGs  are equivalently representations of a chain of binary relations. Every relation of the cobweb poset chain is biunivocally represented by the uniquely designated  \textbf{complete} bipartite digraph-a digraph which is a di-biclique  designated by the very  given sequence. The cobweb poset is then to be identified with a chain of di-bicliques i.e. by definition - a chain of complete bipartite one direction digraphs. Any chain of relations is therefore obtainable from the cobweb poset chain of complete relations  via deleting  arcs (arrows) in di-bicliques.\\ 
Let us underline it again : \textit{any chain of relations is obtainable from the cobweb poset chain of  complete relations via deleting arcs in di-bicliques of the complete relations chain.} For that to see note that any relation  $R_k$ as a subset of  $A_k \times  A_{k+1}$ is represented by a
one-direction bipartite digraph  $D_k$.  A "complete relation"  $C_k$ by definition is identified with its one direction di-biclique graph  $d-B_k$
Any  $R_k$ is a subset of  $C_k$. Correspondingly one direction digraph  $D_k$ is a subgraph of an one direction digraph of  $d-B_k$.\\
The one direction digraph of  $d-B_k$ is called since now on \textbf{the di-biclique }i.e. by definition - a complete bipartite one direction digraph.
Another words: cobweb poset defining di-bicliques are links of a complete relations' chain [6]. 

\vspace{0.2cm}

\noindent Because of that the cobweb posets in the family of all chains of relations unavoidably is of principle 
importance being the most overwhelming case of relations' infinite chains  or their finite parts (i.e. vide - subposets). 

\noindent The intuitively transparent names used above (a chain of di-cliques etc.) are to be the names of  defined 
objects in what follows. These are natural correspondents to their undirected graphs relatives
as we can view a directed graph as an undirected graph with arrowheads added. 

\noindent The purpose of this report is to inform several questions intriguing on their own
apart from being fundamental for a new class of DAGs introduced below. 
Specifically this concerns problems which arise naturally in connection with
a new join combinatorial interpretation of all classical $F-binomial$
coefficients -  Newton binomial, Gaussian $q$-binomial and Fibonomial
coefficients included. This report is based on  [6-7] from which definitions and description
of these new DAG's are quoted for the sake of self consistency and last results are taken from [1]. 
\noindent Applications of new cobweb posets` originated Whitney numbers from [8] 
such as extended Stirling or Bell numbers are expected to be of at least such
a significance in applications to linear algebra of formal series 
as Stirling numbers, Bell numbers or their $q$-extended correspondent already are
in the so called  coherent states physics [9,10] (see [13] for abundant references on this subject and other applications such as in [11,12]).

\vspace{0.2cm}

\noindent \textbf{The problem to be the next.}
\noindent As cobweb subposets $P_n$ are   vertex partitioned into  antichains $\Phi_r$ for $r = 0, 1, ...,n$
which we call levels -  a question of canonical importance arises.
\noindent Let $\left\{P_n\right\}_{n\geq0}$  be the sequence of finite cobweb subposets (see- below). 
What is the form and properties  of $\left\{P_n\right\}_{n\geq0}$'s  characteristic polynomials $\left\{\rho_n(\lambda)\right\}_{n\geq0}$  [15,17]?
\noindent For example - are these related to umbra polynomials? What are recurrence relations defining the $\left\{\rho_n(\lambda)\right\}_{n\geq0}$ family ? This is being now under investigation with a progress to be announced soon by Ewa Krot-Sieniawska [4,5]
- the member of our Rota Polish Seminar Group. 
\noindent The recent papers on DAGs related to this article and its clue references [6,7,8] apart from [1] are [4] and [14].

\vspace{0.1cm}

\noindent \textbf{Computation and Characterizing Problems}.
 
\noindent In the next section we define cobweb posets. Their 
examples  are given [6-8,1]. A join combinatorial interpretation of cobweb
posets` characteristic binomial-like coefficients is provided too following the source papers of the first author.
This simultaneously means that we have join combinatorial interpretation of fibonomial coefficients and 
all incidence coefficients of reduced incidence algebras of full binomial type [16]. 

\noindent In [6,7] the first author  had formulated three problems:  characterization and/or computation of cobweb admissible sequences \textbf{Problem 1}, cobweb layers partition characterization and/or computation \textbf{Problem 2} and the
GCD-morphic sequences characterizations and/or computation \textbf{Problem 3}  all three  interesting on their own.
\noindent Here we report more on one of those problems  [1].  This is a tiling problem. Our information on tiling problem refers to proofs of tiling's existence for some cobweb-admissible sequences as in [1]. There the second author shows that not all cobwebs admit tiling as defined below and 
provides examples  with proofs of families of cobweb posets which  admit tiling.

\section{Cobweb posets as complete bipartite  digraph sequences and combinatorial interpretation [6]}

\noindent \textbf{Cobweb posets as complete bipartite  digraph sequences.} 

\noindent Cobweb posets and their natural subposets $P_n$ are 
graded posets. They are vertex partitioned into  antichains $\Phi_k$ 
for $k = 0, 1, ..., r,...$ (where $r$ is a nonnegative integer) such that for each $\Phi_r$, all of the elements
covering $x$ are in $\Phi_{r+1}$ and all the elements covered by $x$ are in $\Phi_r$. We shall call
the $\Phi_n$ the  $n-th$-level.  $P_n$ is then $(n+1)$ level ranked poset. We are now in a position to 
observe that the cobweb posets may be identified with a chain of di-bicliques i.e. by definition - a chain of complete bipartite one direction digraphs. This is outstanding property as then any chain of relations is obtainable from the corresponding cobweb poset chain of complete relations just by deleting arcs in di-bicliques of this complete relations chain. Indeed.  For any given natural numbers valued sequence the graded (layered) cobweb posets` DAGs  are equivalently representations of a chain of binary relations. Every relation of the cobweb poset chain is bi-univocally represented by the uniquely designated  \textbf{complete} bipartite digraph - a digraph which is a di-biclique  designated by the very  given sequence. The cobweb poset may be therefore identified with a chain of di-bicliques i.e. by definition - a chain of complete bipartite one direction digraphs. Say it again, any chain of relations is obtainable from the cobweb poset chain of complete relations  via deleting  arcs (arrows) in chains` di-bicliques elements.\\ 

\noindent The above intuitively transparent names (a chain of di-cliques etc.) are the names
of the below  defined objects. These objects are natural correspondents of their undirected graphs relatives as
we can view a directed graph as an undirected graph with arrowheads added. Let us start with primary notions remembering that any cobweb subposet $P_k$ is a DAG of course. 

\noindent A bipartite digraph is a digraph whose vertices can be divided into two disjoint sets $V_1$ and $V_2$ such that every arc  connects a vertex in $V_1 $  and a vertex  in $V_2 $. Note that  there is no arc between two vertices in the same independent set  $V_1 $ or $V_2 $. No two nodes of the same partition set are adjacent.

\noindent A one direction bipartite digraph is a bipartite digraph such that every arc  originates at  a node in $V_1 $  and  terminates at a node in $V_2 $. The extension of "being one direction" to $k$-partite digraphs  is automatic. Note that  there is no arc  between two vertices in the same set.
Intuitively  - one may  color the nodes of a bipartite digraph black and blue such that no arc  exists between like colors.

\noindent A  $k$-partite not necessarily one direction  digraph $D$ is obtained from a $k$-partite undirected  $k$-partite 
graph $G$ by replacing every edge  ${xy}$  of  $G$  with the arc $\left\langle xy \right\rangle$, arc $\left\langle yx \right\rangle$  or both  $\left\langle xy \right\rangle$ and $\left\langle yx \right\rangle$. The partite sets of  $D$ are the partite sets of $G$.

\begin{defn}. 
A simple directed graph $G = (V,E)$ is called bipartite if there exists a partition $V = V_1 + V_2$ of the vertex set
$V$  so that every edge (arc) in E is incident with  $v_1$  and  $v_2$ for some $v_1$  in $V_1$  and $v_2$ in $V_2$.
It is complete if  any node from $V_1$  is adjacent to all nodes of  $V_2$.
\end{defn}
We shall denote our special case biparte one direction digraphs as follows $G = (V_1+V_2,E) \equiv  L_{{F_k},F_{k+1}}$ to
inform that the partition has parts $V_1 $ and $V_2 $  where $\left|V_1 \right|= F_k$  and  $\left|V_2\right|= F_{k+1}.$

\noindent \textbf{Notation.} The  special $k$ - partite $\equiv$  $k$ - level one direction  digraphs considered here shall be
denoted  by the corresponding symbol $L_{p,q,...,r}$. The complete $k$ -partite $\equiv$  $k$ -level one direction complete digraph is coded  $K_{p,q,...,r}$  as in the non directed case (with no place for confusion because only one direction  digraphs are to be considered in what follows). 
\noindent $L_{{F_k},F_{k+1},...,F_n}\equiv \langle\Phi_k \rightarrow \Phi_n \rangle $ denotes  $(n-k+1)$ - partite one direction digraph  $\equiv (n-k+1)$ - level one direction digraph $ n\geq k$ whose partition has ( antichains from $\Pi$) the parts $\Phi_k, \Phi_{k+1},...,\Phi_n,\;  \left|\Phi_k \right|= F_k,  \left|\Phi_{k+1} \right|= F_{k+1},…,\left|\Phi_n \right|= F_n.$ Any  cobweb one-layer  $\langle\Phi_k \rightarrow \Phi_{k+1} \rangle,\quad k<n,\quad k,n \in N\cup\{0\}\equiv Z_\geq$  is a complete bipartite one direction digraph  $L_{{F_i},F_{i+1}} = K_{{F_i},F_{i+1}}$  i.e. by definition it is the \textit{di-biclique}. The cobweb one-layer vertices constitute bipartite set $V(\langle\Phi_k \rightarrow \Phi_{k+1} \rangle)=\Phi_k  + \Phi_{k+1}$ while edges are all those arcs incident with the two  antichains nodes in the poset $\Pi$ graph representation.\\
Any  cobweb subposet  $ P_n \equiv\langle\Phi_0 \rightarrow \Phi_n \rangle \equiv L_{F_0,F_1,...,F_n}, n\geq 0$  is a $(n+1)$-partite = $(n+1)$-level one direction digraph. We shall keep on calling  a complete bipartite one direction digraph a  \textbf{di-biclique} because it is a special kind of bipartite one direction digraph, whose every vertex of the first set is connected by an arc originated in this very node to \textbf{every} vertex of the second set of the given bi-partition. 
\noindent Any  cobweb layer $\langle\Phi_k \rightarrow \Phi_n \rangle,\quad k\leq n,\quad k,n \in N\cup\{0\}\equiv Z_\geq$  is a one direction  the $(n-k+1)$-partite $\equiv$ $(n-k+1)$-level one direction DAG  $L_{{F_k},F_{k+1},...,F_n}$ with an additional defining property: it is a chain of di-bicliques  for  $k<n$.  Since now on we shall identify both:  
                       $$\langle\Phi_k \rightarrow \Phi_n \rangle= L_{{F_k},F_{k+1},...,F_n}.$$
Note:  $P_n$ $\equiv$ $\langle\Phi_0 \rightarrow \Phi_n \rangle= L_{{F_0},F_{k+1},...,F_n},\; n\geq 0$  and 
$\langle\Phi_k \rightarrow \Phi_n \rangle \neg = K_{{F_k},F_{k+1},...,F_n}$   for $n>k+1.$\\

\begin{observen} [6] 
$$ \left|E(\langle\Phi_k \rightarrow \Phi_{k+m} \rangle) \right|= {\sum}  ^{m-1}_{i=0}F_{k+i}F_{k+i+1},$$
where  $E(G)$ denotes the set of edges of a graph $G$ (arcs of a digraph $G$).
\end{observen}

\noindent The following  property $(*)$

$$ (*)  \ \ \  \langle\Phi_k \rightarrow \Phi_{k+1}\rangle \equiv L_{{F_k},F_{k+1}}=K_{{F_k},F_{k+1}},\ \ k = 0,1,2,...$$
i.e. $L_{{F_k},F_{k+1}}$  is a di-biclique for $k = 0,1,2,...$ , might be considered the definition of  $F_0$ rooted 
$F$-cobweb graph $\Pi$ or in short $P$ if $F$-sequence has been established. The cobweb poset $\Pi$ is being thus identified with a \textbf{chain of di-bicliques}.  The usual convention is to  choose  $F_0 =1.$  One may relax this constrain, of course.

\noindent Thus any  cobweb one-layer  $\langle\Phi_k \rightarrow \Phi_{k+1}\rangle$ is a complete one direction  bipartite digraph    i.e.  a  di-biclique. This is how the definition of the $F$-cobweb graph $\Pi(F)$ or in short $P$ - has emerged.

\vspace{2mm}

\noindent For infinite cobweb poset $\Pi$  with the set of vertices $P =V(\Pi)$ one has an obvious $\delta(G)=2$ domatic vertex 
partition of this $F$-cobweb poset.

\begin{defn}. 
A subset $D$ of the vertex set $V(G)$ of a digraph $G$ is called dominating in $G$, if each vertex of $G$ either is in $D$, or is adjacent to a vertex of $D$. Adjacent means that  there exists an originating or terminating arc in between the two- any node from $G$  outside $D$  and  a node from $D$.
\end{defn}

\begin{defn}. 
A domatic partition of $V$ is a partition of $V$ into dominating sets, and the number of these dominating sets is called the size of such a partition.  The domatic number  $ \delta(G)$ is the maximum size of a domatic partition.
\end{defn}

\vspace{2mm}

\noindent An infinite cobweb poset $\Pi$  with the set of vertices $P =V(\Pi)$  has the $\delta(\Pi)=2$ domatic vertex 
partition, namely a $mod \ 2$ - partition. 

\noindent $V = V_0 \cup V_1$ where  $V_0 = \bigcup _{k= 2s+1}\Phi_k , \  s = 0,1,2,... $  - ("black levels");
\noindent and  $V_1 = \bigcup _{k= 2s}\Phi_k , \  s = 0,1,2,... $  - ("blue levels").

\noindent Note:  Natural $mod \ n$ partitions of  the cobweb poset`s set of vertices $V(\Pi)=P$  ($n$ colours), 
$P = V_0 \cup V_1 \cup V_2 \cup ... \cup V_{n-1}$ ,  $V_i =\bigcup _{k= 2s+i} , s = 0,1,2,...,i \in Z_n = {0,1,...,n-1}$
for $n>2$ are \textbf{not domatic}.

\noindent  Cobweb one-layer or more than one-layer subposets  $\langle\Phi_k \rightarrow \Phi_{k+m}\rangle \equiv L_{{F_k},F_{k+1}}=K_{{F_k},F_{k+1}},\ \ k = 0,1,2,...$ have also correspondent, obvious the $ \delta(G)=2$ domatic partitions for $m>0.$

\vspace{4mm}

\noindent \textbf{Combinatorial interpretation [7,19,24,25]}. 

\vspace{1mm}

\textbf{2.1. $F-binomial$ coefficients}.
\noindent The source papers are [9-13] from which indispensable definitions and notation are taken for granted including Kwa\'sniewski [9,10] upside - down notation $n_F \equiv F_n$  being used   for dipper than mnemonic reasons - as it  the case with widespread: a]  Gaussian numbers $n_q$ in finite geometries  and the so called "quantum groups" ([8,9]) or b] their $p,q$ cognates $n_{p,q} = \sum_{j=0}^{n-1}{p^{n-j-1}q^j}$ , $n_q  = n_{1,q} $ .

\noindent Given any sequence $\{F_n\}_{n\geq 0}$ of nonzero reals
($F_0 = 0$ being sometimes acceptable as $0! = F_0! = 1.$)
one defines  its corresponding  binomial-like $F-nomial$
coefficients as in Ward`s Calculus of sequences [18] as follows.

\begin{defn}.$(n_{F}\equiv F_{n}\neq 0,\quad n > 0)$
$$
\left( \begin{array}{c} n\\k\end{array}
\right)_{F}=\frac{F_{n}!}{F_{k}!F_{n-k}!}\equiv
\frac{n_{F}^{\underline{k}}}{k_{F}!},\quad \quad n_{F}!\equiv n_{F}(n-1)_{F}(n-2)_{F}(n-3)_{F}\ldots 2_{F}1_{F};$$
$$ 0_{F}!=1;\quad n_{F}^{\underline{k}}=n_{F}(n-1)_{F}\ldots (n-k+1)_{F}. $$
\end{defn}
\noindent We have made above  an analogy driven identifications in the spirit of  Ward`s Calculus of sequences [18]]. 
Identification $n_{F}\equiv F_{n}$ is  the notation used in extended Fibonomial Calculus case [9-13,4] being also there inspiring as $n_F$ mimics $n_q$ established notation for Gaussian integers exploited in much elaborated family of various applications including quantum physics (see [9,10,13] and references therein). 

\vspace{2mm}

\noindent The crucial and elementary observation now is that an eventual
cobweb poset or any combinatorial interpretation of $F$-binomial
coefficients makes sense \textit{ not for arbitrary} $F$
sequences as $F-binomial$ coefficients should be nonnegative
integers (hybrid sets are not considered here).

\vspace{2mm}

\begin{defn}.[7,6,19,24,25]
A natural numbers` valued sequence $F = \{n_F\}_{n\geq 0}$, $F_0 =1$ is called 
cobweb-admissible iff
$$ \left( \begin{array}{c} n\\k\end{array}\right)_{F}\in N_0\quad
for \quad k,n\in N_0.$$
\end{defn}
$F_0 = 0$ being sometimes acceptable as $0_F! \equiv F_0! = 1.$
\vspace{2mm}

\noindent Incidence coefficients of any reduced incidence algebra of full 
binomial type [16] immensely important for computer science are computed exactly with their correspondent cobweb-admissible sequences.
These include binomial (Newton) or  $q$- binomial  (Gauss) coefficients. For other $F$-binomial 
coefficients - computed with cobweb admissible sequences - see in what follows after Observation 3.

\vspace{2mm}

\noindent \textbf{Problem 1}. \textit{ Find effective characterizations  and/or
an algorithm to produce the cobweb admissible sequences i.e. find all examples.} [6,7,19]

\vspace{2mm}
Very recently the second author have proved (a note in preparation) that the following is true.
 
\begin{theoremn}[Dziemia\'nczuk] 
Any cobweb-admissible sequence $F$ is at the point product [1] of primary cobweb-admissible sequences $P(p)$.
\end{theoremn}

\noindent Right from the definition of $P$ via its
Hasse diagram pictures in [6-8] the important observations follow which lead to a specific, new joint
combinatorial interpretation of cobweb poset`s characteristic binomial-like coefficients  [6-8].  

\vspace{3mm}

\begin{observen} [6,7,19,24,25]
The number of maximal chains starting from The Root  (level
$0_F$) to reach any point at the $n-th$ level  with $n_F$ vertices
is equal to $n_{F}!$.
\end{observen}

\begin{observen},  $(k>0)$  [6,7,19,24,25]
The number of all maximal chains in-between $(k+1)-th$ level $\Phi_{k+1}$
and the $n-th$ level $\Phi_n$ with $n_F$ vertices
is equal to $n_{F}^{\underline{m}}$, \quad where $m+k=n.$ 
\end{observen} 

\vspace{1mm}

\noindent Indeed. Denote the number of ways to get along maximal chains from 
\textit{any fixed point} (the leftist for example) in $\Phi_k$ to any vertex 
in  $\Phi_n , n>k$ with the symbol\\
  $$[\Phi_k \rightarrow \Phi_n]$$
then obviously we have ( $[\Phi_n \rightarrow \Phi_n]\equiv 1)$:\\
           $$[\Phi_0 \rightarrow \Phi_n]= n_F!$$ 
and
$$[\Phi_0 \rightarrow \Phi_k]\times [\Phi_k\rightarrow \Phi_n]=
[\Phi_0 \rightarrow \Phi_n].$$

\vspace{2mm}

\noindent For the purpose of a new joint combinatorial interpretation
of $F-sequence-nomial$ coefficients ({\it F-nomial} - in short)
let us consider all finite \textit{"max-disjoint"} sub-posets
rooted at the $k-th$ level at any fixed vertex $\langle
r,k \rangle, 1 \leq r \leq k_F $  and ending  at corresponding
number of vertices at the $n-th$ level ($n=k+m$) where the
\textit{max-disjoint} sub-posets are defined below.

\vspace{2mm}

\begin{defn}. [6,7,19,24,25]
Two posets are said to be max-disjoint if considered as sets of maximal chains 
they are disjoint i.e. they have no maximal chain in common. An equipotent copy of $P_m$ 
[`\textbf{equip-copy}'] is  defined as such a maximal chains family  
equinumerous with $P_m$ set of maximal chains that the it constitutes a sub-poset
with one minimal element.
\end{defn}
We shall proceed with deliberate notation coincidence anticipating coming observation.
\begin{defn}.
Let us denote the number of all mutually max-disjoint equip-copies of $P_m$  rooted
at any  fixed vertex $\langle j,k \rangle , 1\leq j \leq k_F $ of $k-th$ level  with the symbol
$$ \left( \begin{array}{c} n\\k\end{array}\right)_{F}.$$
\end{defn}
One uses here the  customary convention:  $\left(
\begin{array}{c} 0\\0\end{array}\right)_{F}=1$ and $\left(
\begin{array}{c} n\\n\end{array}\right)_{F}=1.$

\vspace{2mm}
\noindent Compare the above with the Definition 4 and the Definition 10.

\vspace{1mm}

\noindent The number of ways to reach an upper level from a
lower one along any of  maximal chains  i.e.  the number of all
maximal chains from the level $\Phi_{k+1}$ to the level $\Phi_n ,\quad k>n$ 
is equal to 
  $$ [\Phi_k \rightarrow \Phi_n]= n_{F}^{\underline{m}}.$$

\noindent Therefore we have

\begin{equation}
\left( \begin{array}{c} n\\k\end{array}\right)_{F} \times [\Phi_0
\rightarrow \Phi_m] = [\Phi_k \rightarrow \Phi_n]=
n_{F}^{\underline{m}}
\end{equation}
where  $[\Phi_0 \rightarrow \Phi_m]= m_F!$ counts the number of
maximal chains in any equip-copy of  $P_m$. With this in mind we see
that the following holds  [6,7,19,24,25].

\vspace{3mm}

\begin{observen} Let  ${n,k}\geq 0$.
Let $n = k+m$. Let $F$ be any cobweb admissible sequence.
Then the number of mutually max-disjoint  equip-copies i.e.  sub-posets
equipotent to $P_{m}$ , rooted at the same \textbf{ fixed} vertex of  $k-th$ level 
and ending at the n-th level is equal to

$$\frac{n_{F}^{\underline{m}}}{m_{F}!} =
\left( \begin{array}{c} n\\m\end{array}\right)_{F}$$
$$ = \left(\begin{array}{c} n\\k\end{array} \right)_{F}=
\frac{n_{F}^{\underline{k}}}{k_{F}!}. $$
\end{observen} 

\noindent The immediate natural question now is

$$\Big\{{{\eta} \atop {\kappa}}\Big\}_{const}= ?$$ 
\vspace{2mm}
\noindent i.e.  the number of partitions with block sizes all equal to const = ?

\vspace{2mm}

\noindent where here $const=\lambda =  m_{F}!$ and 

$$ \eta =  n_{F}^{\underline{m}},\ \  \kappa = \left(\begin{array}{c} n\\k\end{array} \right)_{F} $$
\noindent The    \textit{const} indicates that this is the number of set partitions with  block sizes all equal to 
\textit{const} and we use Knuth notation $\Big\{{{\eta} \atop {\kappa}}\Big\}$ for  Stirling numbers of the second kind.

\vspace{2mm}

\noindent From the formula  (59) in [2] one infers  the Pascal-like matrix answer to the question above.

$$\Big\{{\eta \atop \kappa}\Big\}_{\lambda} = \delta_{\eta,\kappa \lambda}
\frac{\eta !}{\kappa !\lambda !^\kappa}.$$

\noindent  This gives us the rough upper bound for the number of tilings (see [6] for Pascal-like triangles) as we arrive now
to the following intrinsically related problem.

\vspace{2mm}
\noindent \textbf{The partition or tiling Problem 2.}
\noindent Suppose now that  $F$  is a cobweb admissible sequence. Let us introduce  

$$\sigma P_m = C_m[F; \sigma <F_1, F_2,...,F_m>]$$
the equipotent sub-poset obtained from $P_m$ with help of a permutation $\sigma$
of the sequence $<F_1, F_2,...,F_m>$ encoding  $m$ layers
of $P_m$ thus obtaining the equinumerous sub-poset $\sigma P_m$ 
with the sequence $\sigma <F_1, F_2,...,F_m>$ encoding  now $m$ layers of $\sigma P_m$.
Then $P_m = C_m[F; <F_1, F_2,...,F_m>].$ Consider the layer 
$\langle\Phi_k \rightarrow \Phi_n \rangle,\quad k<n,\quad
k,n \in N$ partition into the equal size blocks which are here max-disjoint equi-copies of $P_m, m=n-k+1$.
The question then arises whether and under which conditions the layer may be
partitioned with help of max-disjoint  blocks of the form $\sigma P_m$. And how to visualize this phenomenon?
It seems to be the question of computer art, too. At first - we already know that an answer to 
the main question of such tilings existence - for some sequences $F$ -is in affirmative. 
Whether is it so for all cobweb admissible sequences -we do not know by now.
Some computer experiments done by student Maciej Dziemia\'nczuk [1] are encouraging. 
More than that. The second author in [1] proves tiling's existence for some cobweb-admissible sequences including natural and Fibonacci numbers sequences. He shows also that not all $F$ - designated cobweb posets do admit tiling as defined above. 
However problems:  "how many?" is opened. Let us recapitulate and report on results obtained in [1].

\section {Cobweb posets tiling problem }

\begin{figure}[ht]
\begin{center}
	\includegraphics[width=100mm]{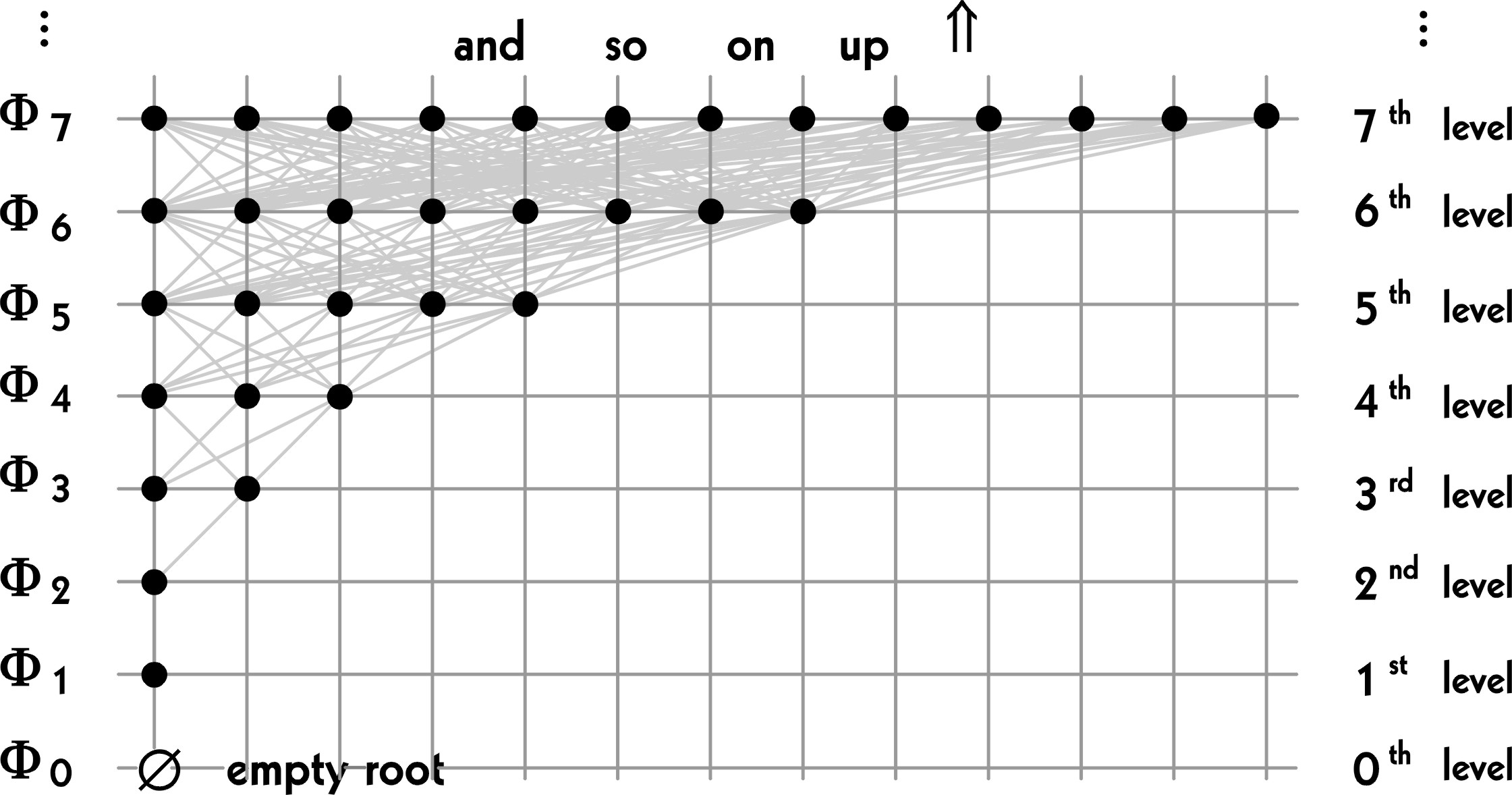}
	\caption{The $s$-th level in $\mathbb{N}\times \mathbb{N}\cup\{0\}$ \label{fig:akkfig0}}
\end{center}
\end{figure}

Let us recall that cobweb poset in its original form [6,7] is defined as a partially ordered graded infinite poset $\Pi = \langle P,\leq\rangle$, designated uniquely by any sequence of nonnegative integers $F = \{ n_F \}_{n \geq 0 }$ and it is represented as a directed acyclic graph (DAG) in the graphical display of its Hasse diagram. $P$ in $\langle P,\leq \rangle$ stays for   set of vertices while $\leq$ denotes partially ordered relation. See Fig. \ref{fig:akkfig0} and note  (quotation from [7,6]):

\begin{quote}
One  refers to $\Phi_s$ as to  the set of vertices at the $s$-th level. The population of the $k$-th level ("\emph{generation}") counts  $k_F$ different member vertices for $k>0$ and one for $k=0$.
Here down (Fig. \ref{fig:akkfig0}) a disposal of vertices on $\Phi_k$  levels is visualized for the case of Fibonacci sequence. $F_0 = 0$ corresponds to the empty root $\{\emptyset\}$.
\end{quote}

In Kwa\'sniewski's cobweb posets' tiling problem one considers finite cobweb sub-posets for which we have finite number of levels in layer $\langle \Phi_k \rightarrow \Phi_n \rangle$, where  $k\leq n$,  $k, n\in\mathbb{N}\cup \{0\} $ with exactly $k_j$  vertices on $\Phi_j$ level $k\leq j \leq n$. For $k=0$ the sub-posets $\langle \Phi_0 \rightarrow \Phi_n \rangle$ are named \emph{ prime cobweb posets} and these are those to be used - up to permutation of levels equivalence - as a block to partition finite cobweb sub-poset.

For the sake of combinatorial interpretation a natural numbers' valued sequence $F$ which determines its' cobweb poset has to be \emph{cobweb-admissible}.

\begin{figure}[ht]
\begin{center}
	\includegraphics[width=75mm]{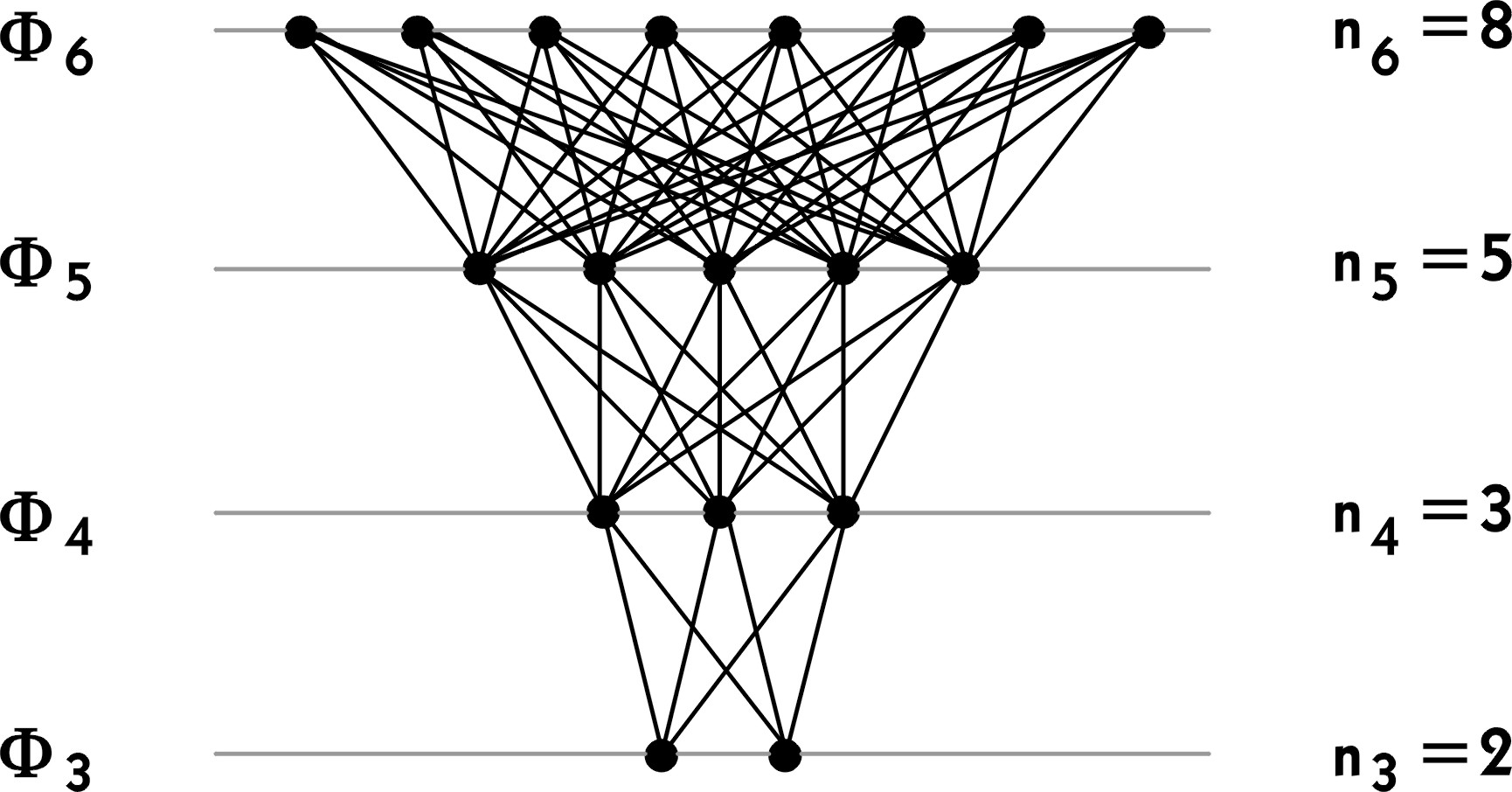}
	\caption{Display of four levels of Fibonacci numbers' finite Cobweb sub-poset}
\end{center}
\end{figure}

\begin{figure}[ht]
\begin{center}
	\includegraphics[width=75mm]{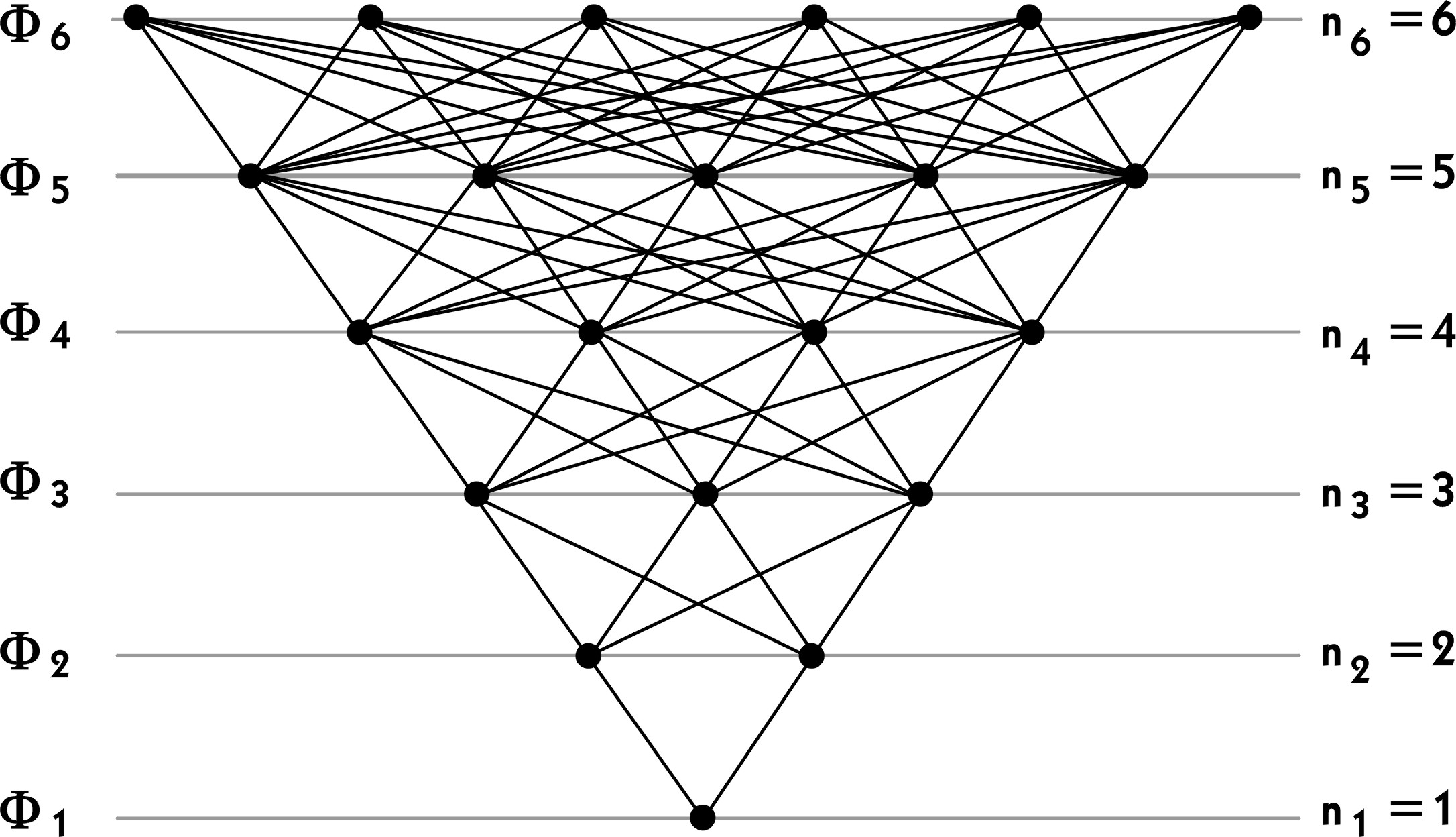}
	\caption{Display of Natural numbers' finite prime Cobweb poset}
\end{center}
\end{figure}

\noindent
$F_0 = 0$ being acceptable as $0_F! \equiv F_0! = 1$. 
We adopt then the convention to call the root $\{\emptyset\}$ the "empty root".

\vspace{0.4cm}

\noindent One of the problems posed in [6-8] is the one, which is the subject of [1].

\vspace{0.4cm}

\noindent \textbf{The tiling problem}

\vspace{0.2cm}

\noindent Suppose now that $F$ is a cobweb admissible sequence. Under which conditions any layer $\langle\Phi_n\!\rightarrow\!\Phi_k\rangle$ may be partitioned with help of max-disjoint blocks of established type $\sigma P_m$? Find effective characterizations and/or find an algorithm to produce these partitions.

\vspace{0.2cm}

The above Kwa\'sniewski [7,6] tiling problem  is first of all the problem of existence of a partition of layer $\langle\Phi_k \rightarrow\Phi_n \rangle$  with max-disjoint blocks of the form $\sigma P_m$  defined as follows:

\begin{displaymath}
	\sigma P_m = C_m [F, \sigma \langle F_1, F_2, \ldots, F_m \rangle ]
\end{displaymath}

\begin{figure}[ht]
\begin{center}
	\includegraphics[width=100mm]{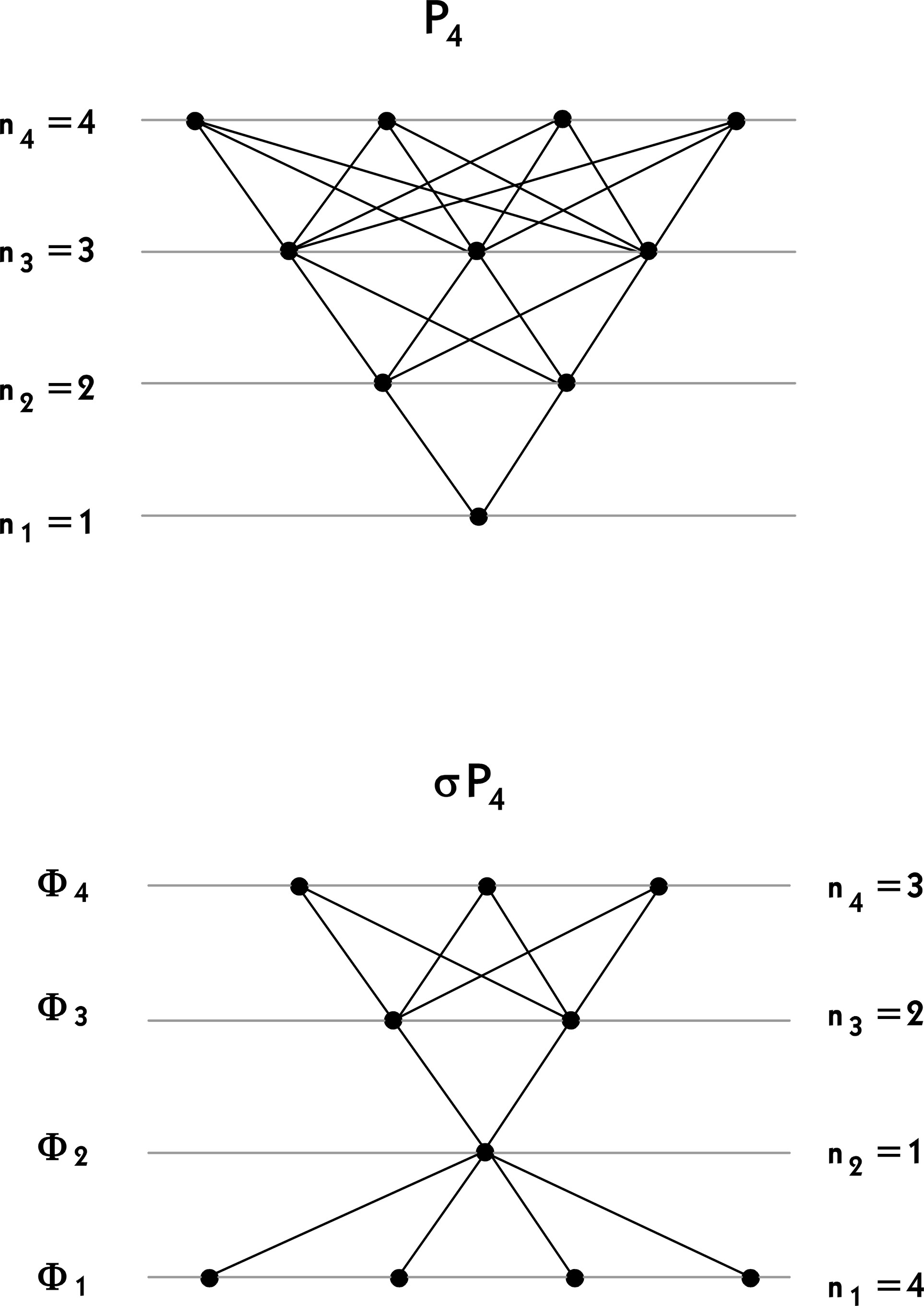}
	\caption{Display of block $\sigma P_m$ obtained from $P_m$ and permutation $\sigma$ \label{fig:permut}}
\end{center}
\end{figure}

It means that the partition may contain only primary cobweb sub-posets or these obtained from primary cobweb poset $P_m$ via permuting its levels as illustrated below (Fig. \ref{fig:permut}).

The second author presents in [1] an algorithm to create a partition of any layer $\langle\Phi_k \rightarrow\Phi_n \rangle$, $k\leq\nolinebreak n$, $k,n \in \mathbb{N}\cup \{0\}$   of finite cobweb sub-poset specified by such  $F$-sequences as Natural numbers and Fibonacci numbers.  In [1] the following  Theorem 1  and Theorem 2 are proved. 

\begin{theoremn}[Natural numbers] 
Consider any layer $\langle \Phi_{k+1} \rightarrow \Phi_n \rangle$ with $m$ levels where $m=n-k$, $k\leq{n}$ and  $k,n\in\mathbb{N} \cup \{0\}$ in a finite cobweb sub-poset, defined by the sequence of \textbf{natural numbers} i.e. $F \equiv \{ n_F \}_{n\geq 0},\ n_F = n,\ \linebreak n\in\mathbb{N} \cup \{0\}$. Then there exists at least one way to partition this layer with help of max-disjoint blocks of the form $\sigma P_m$.
\end{theoremn}

\noindent Max-disjoint means that the two blocks have no maximal chain in common.

\noindent Before proving let us notice that for any $m, k \in \mathbb{N}$ such that $m+k=n$:

\begin{equation}
	n_F = m_F + k_F	\label{eq:1}
\end{equation}

\noindent where $1_F = 1$.

\begin{figure}[ht]
\begin{center}
	\includegraphics[width=60mm]{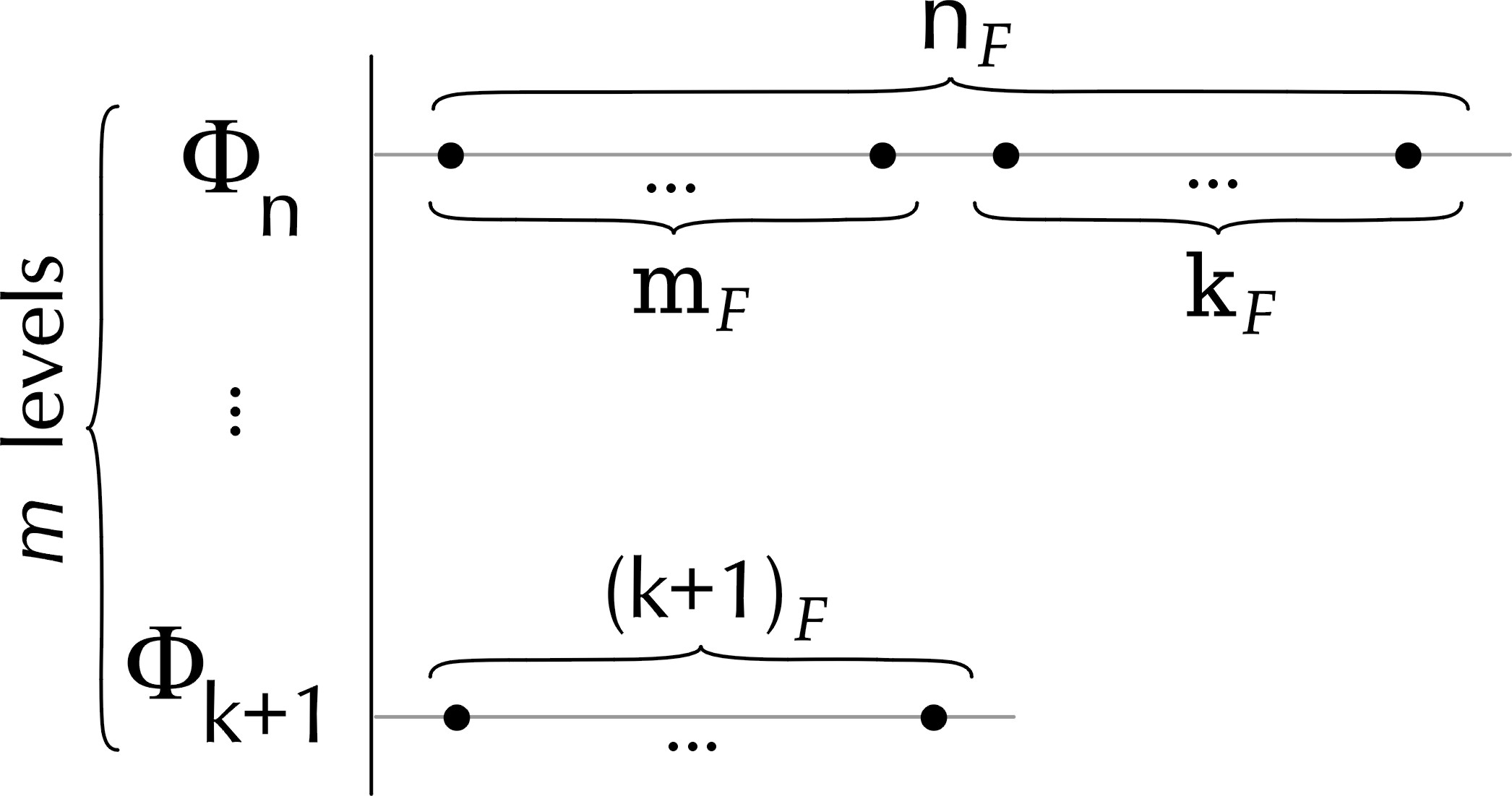}
	\caption{Picture of $m$ levels of Cobweb poset' Hasse diagram \label{fig:steep1} }
\end{center}
\end{figure}

% % PROOF 
\noindent \textbf{P}\textbf{\footnotesize{ROOF}} \textbf{(cprta1) algorithm}

\vspace{0.2cm}

\noindent \textbf{Steep 1}. There are $n_F = m_F + k_F$ vertices on the $\Phi_n$ level. Let us separate them by cutting into two disjoint subsets as illustrated by the Fig.\ref{fig:steep1} and cope at first with $m_F$ vertices (Steep~2).  Then we shall cope with those $k_F$ vertices left (Steep~3).

\begin{figure}[ht]
\begin{center}
	\includegraphics[width=50mm]{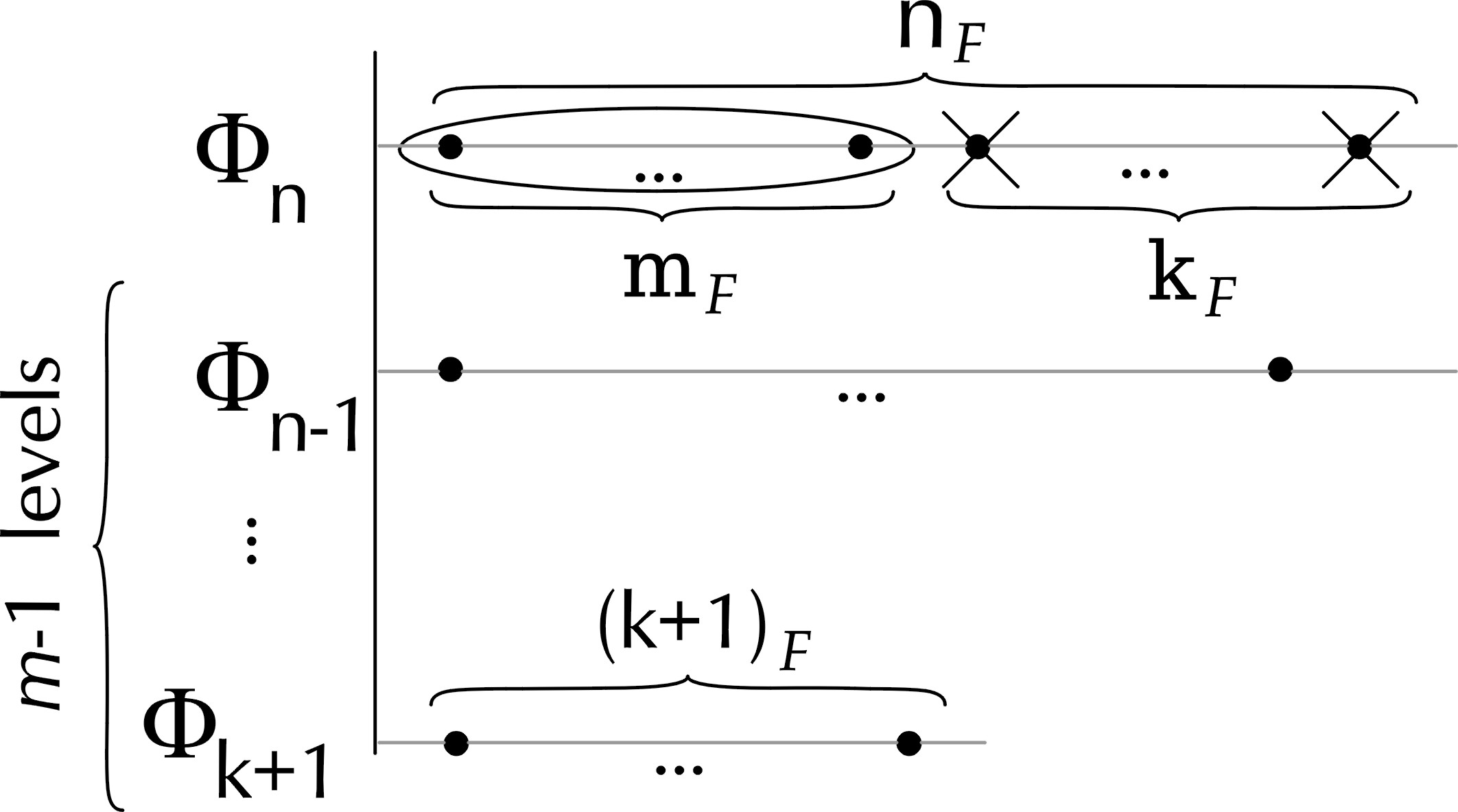}
	\caption{Picture of Steep 2 \label{fig:steep2} }
\end{center}
\end{figure}

\noindent \textbf{Steep 2}. Temporarily we have $m_F$ fixed vertices on $\Phi_n$ level to consider. Let us cover them by $m$-th level of block $P_m$, which has exactly $m_F$ vertices-leafs. What was left is the layer $\langle \Phi_{k+1} \rightarrow \Phi_{n-1} \rangle$ and we might eventually partition it with smaller max-disjoint blocks $\sigma P_{m-1}$, but we need not to do that.  See the next step.

\vspace{0.2cm}
\noindent \textbf{Steep 3}. Consider now the second complementary situation, where we have $k_F$ vertices on $\Phi_n$ level being fixed.  Observe that if we \emph{move} this level lower than $\Phi_{k+1}$ level, we obtain exactly $\langle \Phi_{k} \rightarrow \Phi_{n-1} \rangle$ layer to be partitioned with max-disjoint blocks of the form $\sigma P_m$.  This "\emph{move}" operation is just permutation of levels' order.

\vspace{0.2cm}

The layer $\langle\Phi_{k+1}\!\!\rightarrow\!\Phi_{n}\rangle$ may be partitioned with $\sigma P_m$ blocks if $\langle\Phi_{k+1}\!\rightarrow\nolinebreak\Phi_{n-1}\rangle$ may be partitioned with $\sigma P_{m-1}$ blocks and $\langle\Phi_{k}\!\!\rightarrow \Phi_{n-1} \rangle$ by $\sigma P_m$ again. Continuing these steeps by induction, we are left to prove that $\langle\!\Phi_{k}\!\rightarrow \Phi_{k} \rangle$ may be partitioned by $\sigma P_1$  blocks and $\langle \Phi_{1}\!\!\rightarrow\!\Phi_{m} \rangle$ by $\sigma P_m$ blocks which is obvious $\blacksquare$

\begin{figure}[ht]
\begin{center}
	\includegraphics[width=85mm]{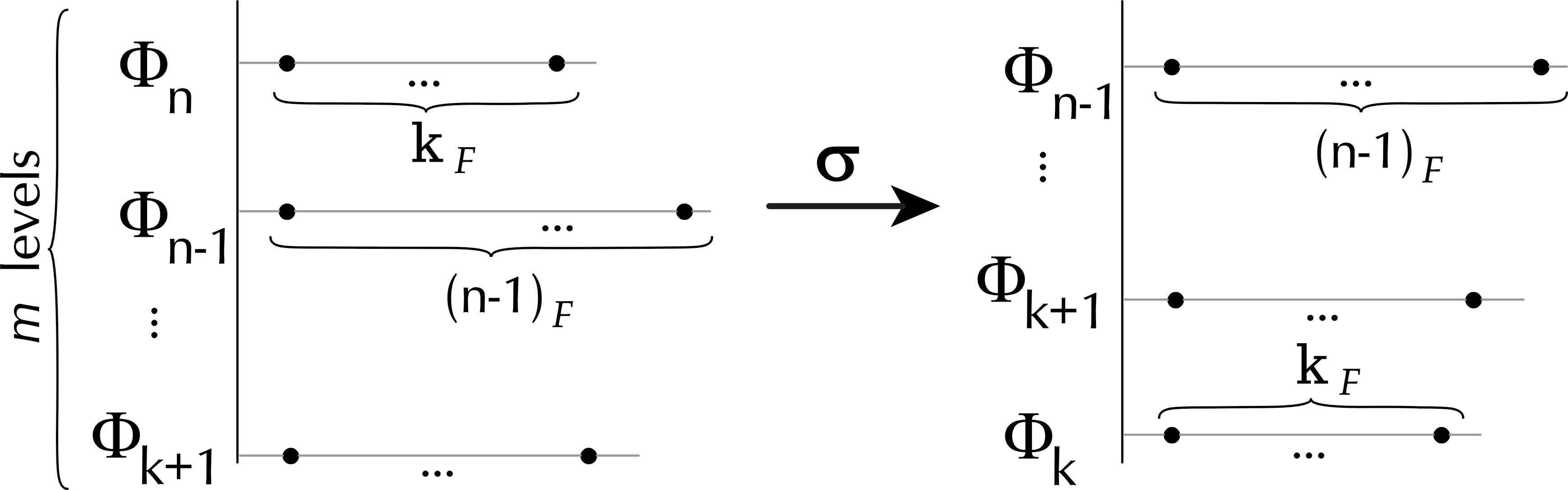}
	\caption{Picture of Steep 3 \label{fig:steep3} }
\end{center}
\end{figure}

\begin{observen}
$\\ $\emph{ We already know from  [7,6]  that the number of max-disjoint equip-copies of $\sigma P_m$, rooted at the same fixed vertex of $k$-th level and ending at the $n$-th level is equal to }
\end{observen}

\begin{displaymath}
	{n \choose k}_F = {n \choose m}_F
\end{displaymath}

\noindent If we cut-separate family of leafs of the layer $\langle\Phi_{k+1}\!\!\rightarrow\Phi_{n}\rangle$, as in the proof of the Theorem 1 then the number of max-disjoint equip copies of $P_{m-1}$ from the Steep 2 is equal to

\begin{displaymath}
	{n-1 \choose k}_F
\end{displaymath}

\noindent However the number of max-disjoint equip copies of $P_m$ from the Steep 3 is equal to

\begin{displaymath}
	{n-1 \choose k-1}_F
\end{displaymath}

\noindent It results in  formula of Newton's symbol recurrence:

\begin{displaymath}
	{n \choose k}_F = {n-1 \choose k}_F + {n-1 \choose k-1}_F
\end{displaymath}

\vspace{0.2cm}
\noindent in accordance with  what was expected for the case  $F = \mathbb{N}$  thus illustrating  the combinatorial interpretation from  [7,6]
in this particular case.  

\vspace{0.4cm}
In the next we adapt Knuth notation for "$F$-Stirling numbers" of the second kind $\Big\{ {n \atop k } \Big\}_F $ as in [6] and also in conformity with Kwa\'sniewski notation for \linebreak $F$-nomial coefficients [9-13,4].
The number of those partitions which are obtained via (cprta1) algorithm shall be denoted by the symbol $\Big\{ {n \atop k } \Big\}_F^1$.

\begin{observen} 
$\\ $\emph{ Let $F$ be a sequence matching (\ref{eq:1}). Then the number $\Big\{ {n \atop k } \Big\}_F^1 $ of different partitions of the layer $\langle\Phi_{k}\!\!\rightarrow\Phi_{n}\rangle$ where $n,k \in \mathbb{N},\ n,k \geq 1$ is equal to:}
\end{observen}

\begin{tabular}{ p{95mm}  p{10mm}}
	\begin{displaymath} 
	\bigg\{ {n \atop k } \bigg\}_F^1  =  {n_F \choose m_F}  \cdot  \bigg\{ {n-1 \atop k } \bigg\}_F^1  	\cdot  \bigg\{ {n-1 \atop k-1 } \bigg\}_F^1
	\end{displaymath}
& 
	\vspace{0.3cm}
	\begin{center} 
		$(S_N)$
	\end{center}
\end{tabular}

\noindent where $\Big\{{n\atop n}\Big\}_F^1 = \Big\{{n\atop n}\Big\}_F = 1$, 
$\Big\{{n\atop 1}\Big\}_F^1 = \Big\{{n\atop 1}\Big\}_F = 1$,
$m = n-k+1$.

\vspace{0.4cm}

% % PROOF
\noindent \textbf{P}\textbf{\footnotesize{ROOF}}

\vspace{0.2cm}
\noindent According to the  Steep 1 of the proof of Theorem 1 we may choose on $\Phi_n$ level $m_F$ vertices out of $n_F$ ones in ${n_F \choose m_F}$ ways.  Next recurrent steps of the proof of  Theorem 1 result in formula ($S_N$) via product rule of counting. $\blacksquare$

\vspace{0.2cm}
\noindent \textbf{Note.} $\Big\{{n\atop k}\Big\}_F^1$ is not the number of all different partitions of the layer $\nolinebreak{\langle\Phi_{k}\!\!\rightarrow\!\Phi_{n}\rangle}$ i.e. $\Big\{{n\atop k}\Big\}_F \geq \Big\{{n\atop k}\Big\}_F^1 $  as computer experiments \cite{6} show. There are much more other tilings with blocks $\sigma P_m$.

\begin{figure}[ht]
\begin{center}
	\includegraphics[width=55mm]{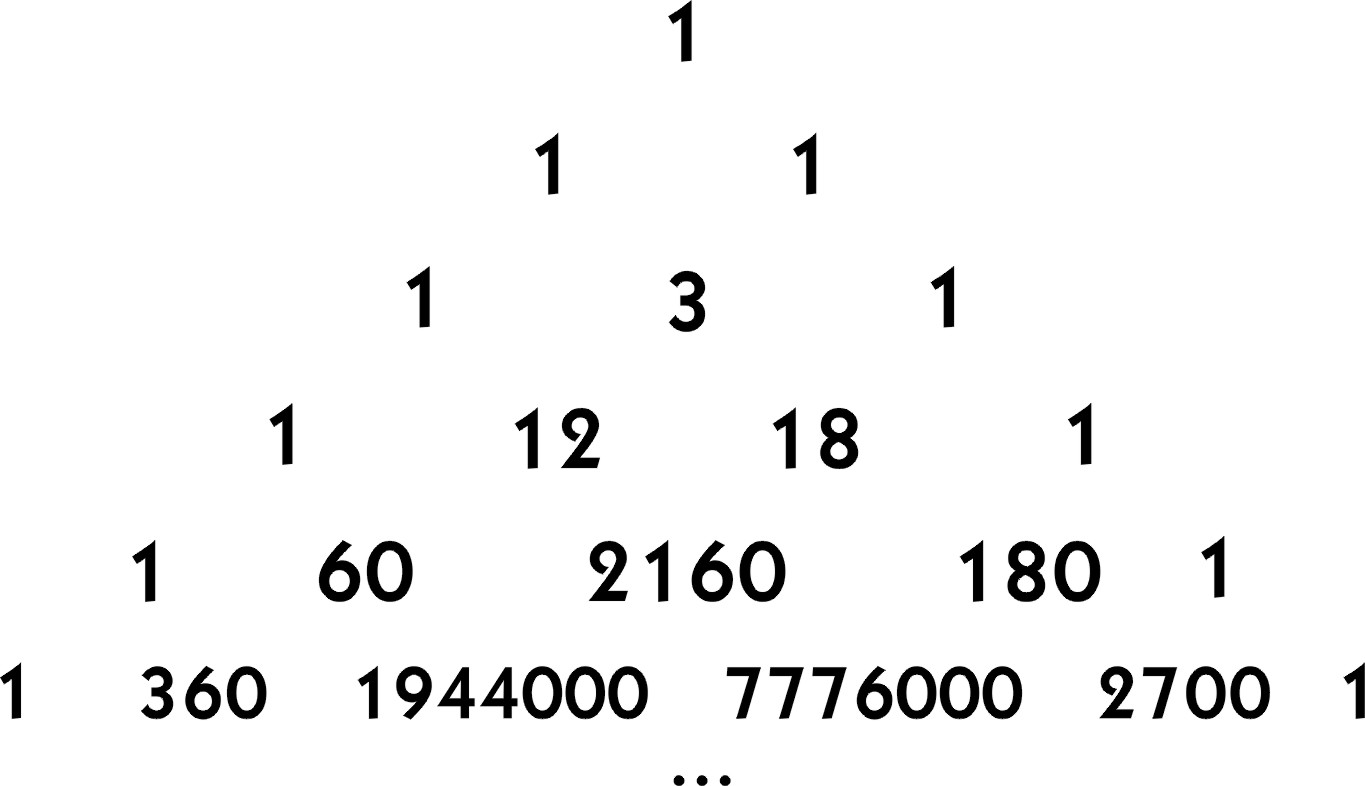}
	\caption{Natural numbers' Cobweb poset tiling triangle of $\Big\{{n\atop k}\Big\}_F^1$ \label{fig:triangle1} }
\end{center}
\end{figure}

\begin{figure}[ht]
\begin{center}
	\includegraphics[width=105mm]{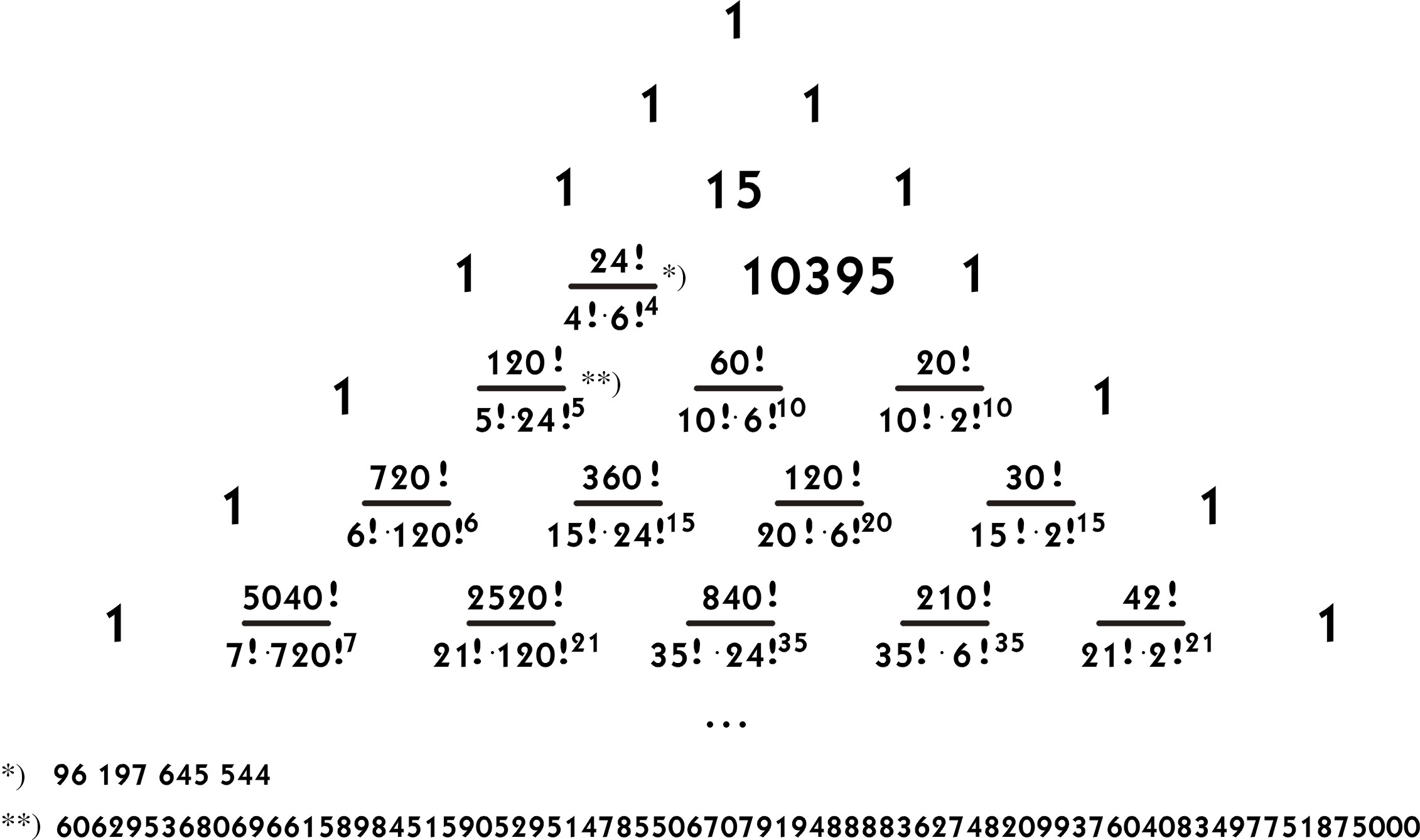}
	\caption{Kwa\'sniewski Natural numbers' cobweb poset tiling triangle of $\Big\{{\eta \atop \kappa}\Big\}_{\lambda}$ \label{fig:akk_triangle1} }
\end{center}
\end{figure}

This is to be  compared  with  Kwa\'sniewski cobweb   triangle   [6] (Fig. \ref{fig:akk_triangle1}) for the infinite triangle matrix  elements

$$\Big\{{\eta \atop \kappa}\Big\}_{\lambda} = \delta_{\eta,\kappa \lambda}
\frac{\eta !}{\kappa !\lambda !^\kappa}$$
\noindent counting the number of partitions with block sizes all equal to $\lambda$.

\noindent Here $const=\lambda =  m_{F}!, m=n-k+1$ and 

$$ \eta =  n_{F}^{\underline{m}},\ \  \kappa = {n \choose k-1}_{F} $$

\noindent The numbers appearing above in $n$-th row, $n>3$ are GIANT numbers as seen from Fig.\ref{fig:akk_triangle1}.

\noindent The inequality 
$\Big\{{n\atop k}\Big\}_F^1 \leq \Big\{{\eta \atop \kappa}\Big\}_{\lambda}$
gives us the rough upper bound for the number of tilings 
with blocks of established type  $\sigma P_m$.

\vspace{0.4cm}
\begin{theoremn}[Fibonacci numbers]
Consider any layer $\langle\Phi_{k+1}\rightarrow\Phi_n \rangle$ with $m$ levels where $m=n-k$, $k\leq n$  and $k,n\in \mathbb{N}\cup\{0\}$ in a finite cobweb sub-poset, defined by the sequence of Fibonacci numbers i.e. $F \equiv \{ n_F \}_{n\geq 0}, n_F \in \mathbb{N} \cup \{0\}$. Then there exists at least one way to partition this layer with help of max-disjoint blocks of the form $\sigma P_m$.
\end{theoremn}

The proof of the Theorem 2 for the Fibonacci sequence $F$ is similar to the proof of Theorem 1. We only need to notice that for any  $m,k \in \mathbb{N}$, $m>1$, $m+k=n$ the following identity takes place:

\begin{equation}
	n_F = (m+k)_F = (k+1)_F\cdot m_F + (m-1)_F\cdot k_F \label{eq:2}
\end{equation}
\vspace{0.4cm}
\noindent where $1_F = 2_F = 1$.

\begin{figure}[ht]
\begin{center}
	\includegraphics[width=60mm]{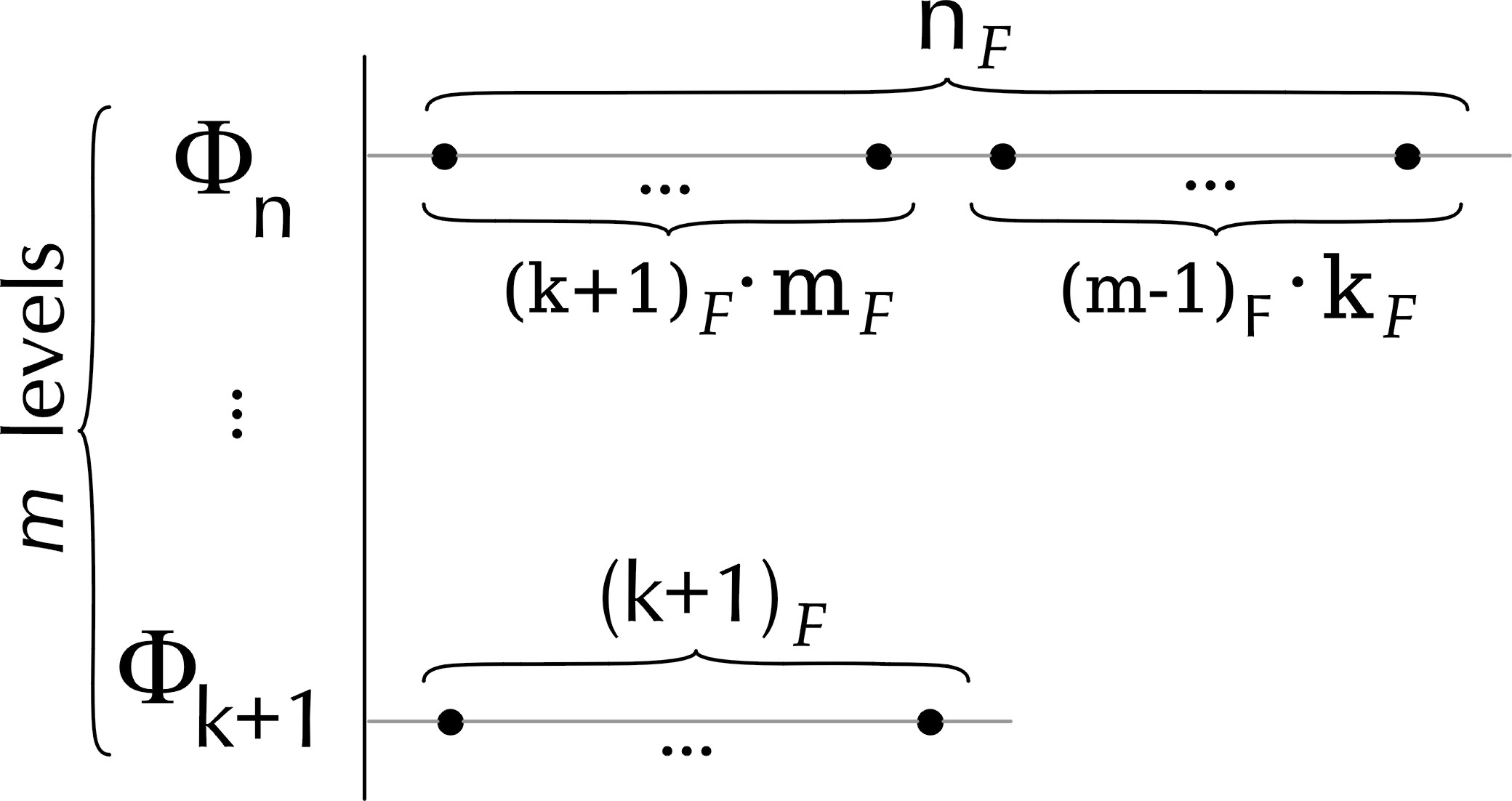}
	\caption{Picture of $m$ levels' layer of Fibonacci Cobweb graph \label{fig:steep1Fib} }
\end{center}
\end{figure}

% % % PROOF
\noindent \textbf{P}\textbf{\footnotesize{ROOF}}

\vspace{0.2cm}

\noindent The number of leafs on the $\Phi_n$ layer is the sum of two summands $\kappa\cdot m_F$  and $\mu\cdot k_F$, where $\kappa=(k+1)_F$, $\mu=(m-1)_F$, (Fig. \ref{fig:steep1Fib}) therefore as in the proof of the Theorem 1 we consider two parts.  At first we have to partition $\kappa$ layers $\langle\Phi_{k+1}\rightarrow\Phi_{n-1}\rangle$ with blocks $\sigma P_{m-1}$ and $\mu$ layers $\langle\Phi_{k}\rightarrow\Phi_{n-1}\rangle$ with $\sigma P_m$. The rest of the proof goes similar as in the case of the Theorem 1 $\blacksquare$

\vspace{0.4cm}
\noindent Theorem 2 is a generalization of Theorem 1 corresponding to $const=\kappa,\mu=1$ case.

\begin{observen}
$\\ $\emph{
The number of max-disjoint equip copies of $P_{m-1}$ which partition $\kappa$ layers $\langle\Phi_{k+1}\rightarrow\Phi_{n-1}\rangle$ is equal to}
\end{observen}

\begin{displaymath}
	\kappa{n-1 \choose k}_F = (k+1)_F{n-1 \choose k}_F
\end{displaymath}
 
However this number of max-disjoint equip copies of $P_m$ which partition $\mu$ layers $\langle\Phi_{k}\rightarrow\Phi_{n-1}\rangle$ is equal to

\begin{displaymath}
	\mu{n-1 \choose k-1}_F = (m-1)_F{n-1 \choose k-1}_F
\end{displaymath}

Therefore the  sum corresponding to the Step 2 and to the Step 3 is the well known recurrence relation for Fibonomial coefficients [11,7,6,4]

\begin{displaymath}
	{n \choose k}_F = (k+1)_F{n-1 \choose k}_F + (m-1)_F{n-1 \choose k-1}_F
\end{displaymath}

\noindent in accordance with what was expected for the case $F$ being now Fibonacci sequence
thus illustrating the combinatorial interpretation from [6,7] in this particular case.

\begin{observen}
$\\ $\emph{
Let $F$ be a sequence matching (\ref{eq:2}). Then the number  $\Big\{ {n \atop k } \Big\}_F^1$ of different partitions of the layer $\langle\Phi_{k}\rightarrow\Phi_{n}\rangle$ where $n,k\in\mathbb{N}, n,k\geq 1$ is equal to:}
\end{observen}

\begin{tabular}{ p{95mm}  p{1cm}}
	\begin{displaymath} 
	\bigg\{ {n \atop {k} } \bigg\}_F^1  =  \frac{F_n!}{(F_m!)^\kappa \cdot (F_{k-1}!)^\mu}
  \cdot  \bigg\{ {n-1 \atop k } \bigg\}_F^1  	\cdot  \bigg\{ {n-1 \atop k-1 } \bigg\}_F^1
	\end{displaymath}
& 
	\vspace{0.3cm}
	\begin{center} 
		$(S_F)$
	\end{center}
\end{tabular}

\noindent where
$\Big\{{n\atop n}\Big\}_F^1 = \Big\{{n\atop n}\Big\}_F = 1$, 
$\Big\{{n\atop n-1}\Big\}_F^1 = \Big\{{n\atop n-1}\Big\}_F = 1$,
$\Big\{{n\atop 1}\Big\}_F^1 = \Big\{{n\atop 1}\Big\}_F = 1$,
$\kappa=k_F, \mu=(m-1)_F, m=n-k+1$, 
$F_n! = 1\cdot 2\cdot \ldots \cdot (n_F-1)\cdot n_F$.

% % % PROOF
\vspace{0.4cm}
\noindent \textbf{P}\textbf{\footnotesize{ROOF}}

\vspace{0.2cm}
\noindent According to the Steep 1 of the proof of Theorem 2 we may choose on  $n$-th level $m_F$ vertices $\kappa$  times and next $(k-1)_F$ vertices $\mu$ times out of $n_F$ ones  in $\frac{F_n!}{(F_m!)^\kappa \cdot (F_{k-1}!)^\mu}$ ways.  Next recurrent steps of the proof of Theorem 2 result in formula  ($S_F$)  via product rule of counting $\blacksquare$

\vspace{0.4cm}
\noindent Observation 4 becomes Observation 2 once we put $const=\kappa,\mu=1$. 

\vspace{0.4cm}
\noindent {\bf Easy example} 

\begin{figure}[ht]
\begin{center}
	\includegraphics[width=70mm]{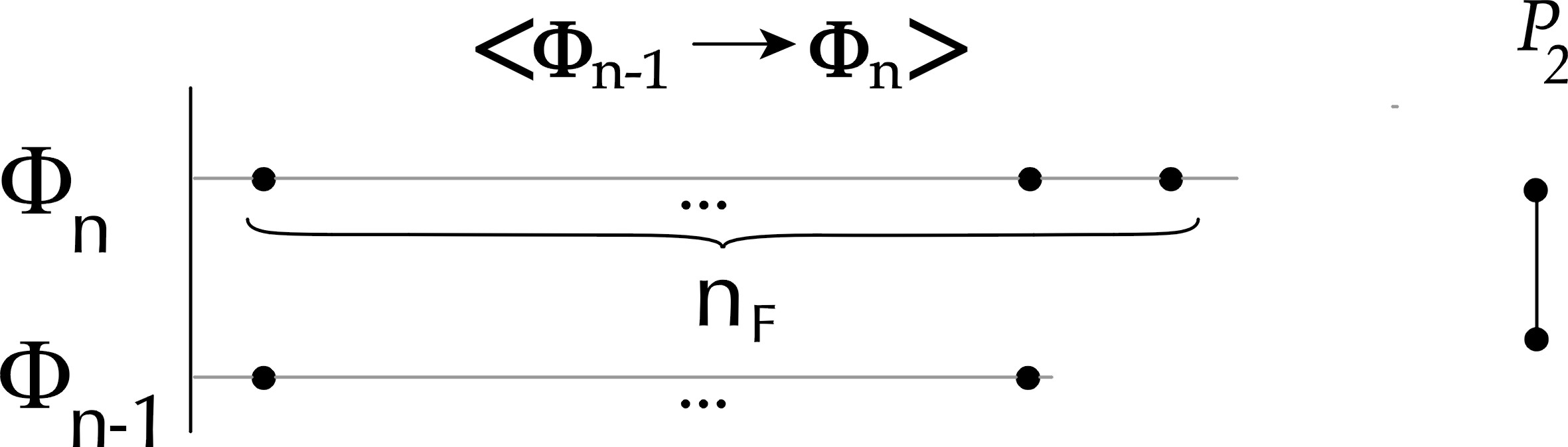}
	\caption{Easy example picture \label{fig:note11}}
\end{center}
\end{figure}

\noindent For cobweb-admissible sequences $F$ such that $1_F = 2_F = 1$,
$\Big\{\!{n\atop n-1}\!\Big\}_F^1\!\!\!=\!\!\Big\{\!{n\atop n-1}\!\Big\}_F\!\!\!=\!1$ as obviously we deal with the perfect matching of the bipartite graph which is very exceptional case (Fig. \ref{fig:note11}).

\noindent \textbf{Note.} As in the case of Natural numbers for $F$-Fibonacci numbers $\Big\{{n\atop 1}\Big\}_F^1$ is not the number of all different partitions of the layer $\langle\Phi_{k}\rightarrow\Phi_{n}\rangle$ i.e. $\Big\{{n\atop k}\Big\}_F \geq \Big\{{n\atop k}\Big\}_F^1$ as computer experiments \cite{6} show. There are much more other tilings with blocks $\sigma P_m$.

\begin{figure}[ht]
\begin{center}
	\includegraphics[width=75mm]{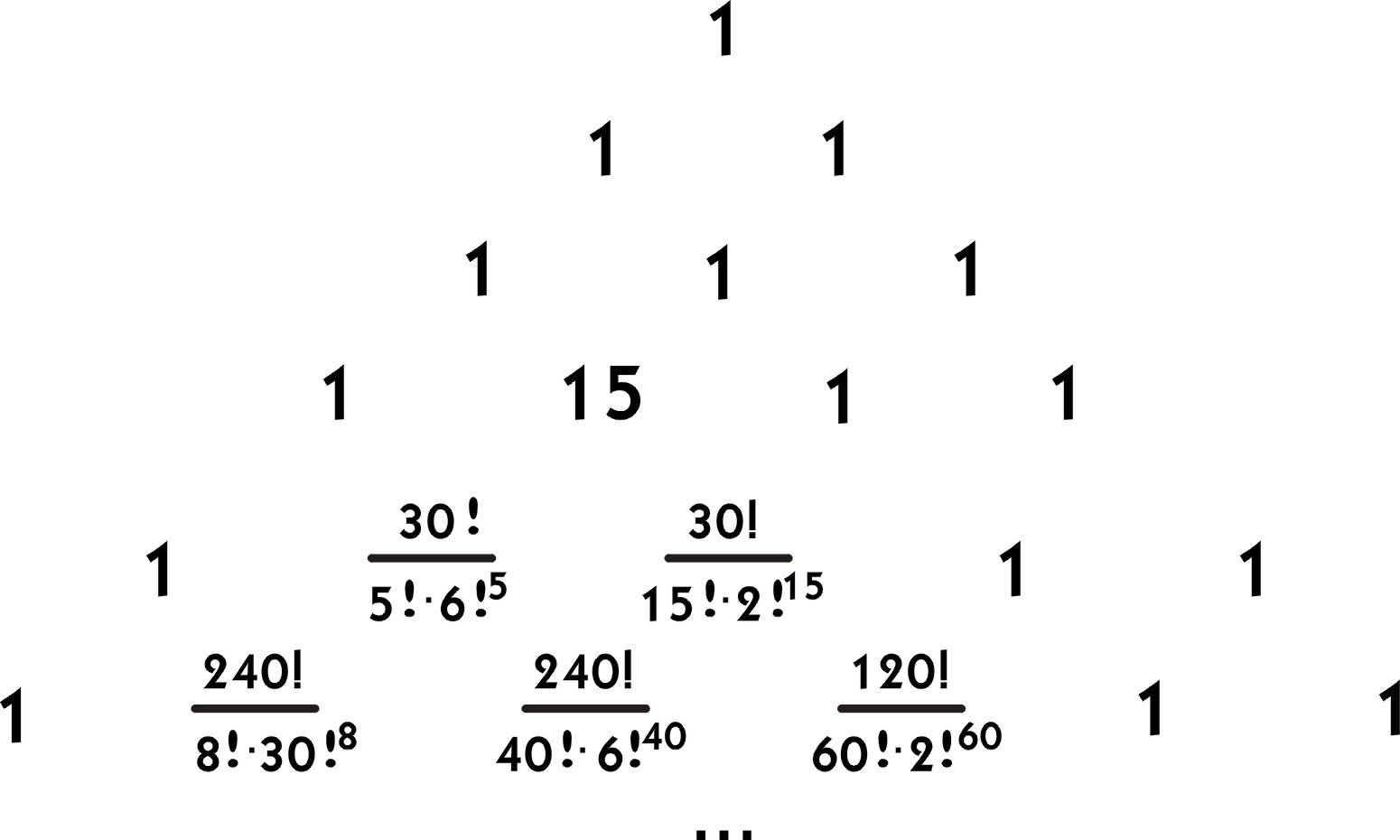}
	\caption{Kwa\'sniewski Fibonacci numbers' cobweb tiling triangle of $\Big\{{\eta \atop \kappa}\Big\}_{\lambda}$ \label{fig:akk_triangle2} }
\end{center}
\end{figure}

This is to be compared with Kwa\'sniewski [6] cobweb triangle for the infinite triangle matrix elements (Fig. \ref{fig:akk_triangle2})

\begin{figure}[ht]
\begin{center}
	\includegraphics[width=60mm]{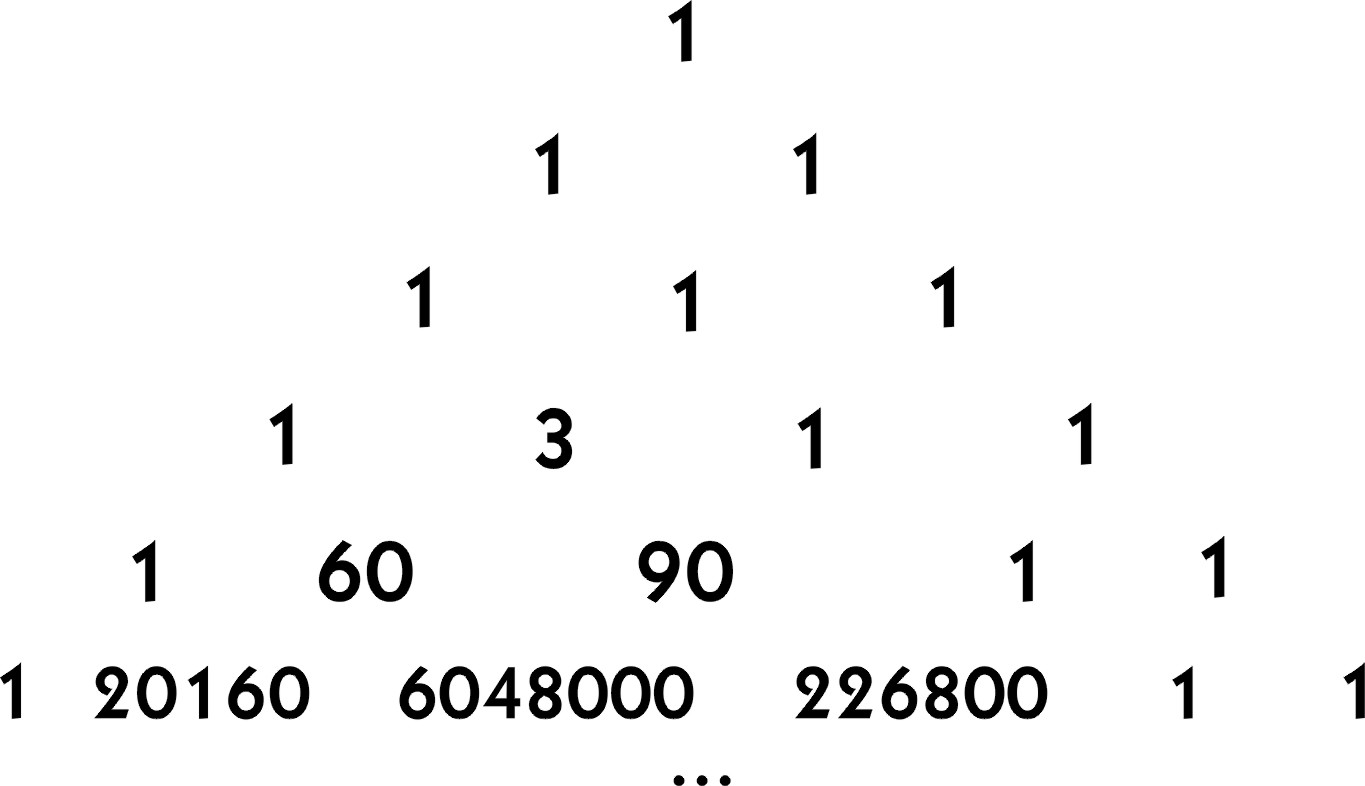}
	\caption{Fibonacci numbers' cobweb poset tiling triangle of $\Big\{{n\atop k}\Big\}_F^1$ \label{fig:triangle1_fib} }
\end{center}
\end{figure}

% % % % % % % % Section
\section {Other tiling sequences}

\begin{defn}
The cobweb admissible sequences that designate cobweb posets with tiling are called cobweb tiling sequences.
\end{defn}

\subsection{Easy examples}
The above method applied to prove tiling existence for Natural and Fibonacci numbers relies on the assumptions (\ref{eq:1}) or (\ref{eq:2}). Obviously these are not the only sequences that do satisfy recurrences (\ref{eq:1}) or (\ref{eq:2}). There exist also other cobweb tiling sequences beyond the above ones with different initial values.
 
There exist also cobweb admissible sequences determining cobweb poset with \emph{no} tiling of the type considered in this paper.

\begin{example}
$n_F = (m+k)_F = m_F + k_F$, $n\geq 1$ \emph{( $0_F$ = corresponds to one "\emph{empty root}" $\{\emptyset\}$)}
\end{example} 

This might be considered a sample example illustrating the method. For example if we choose $1_F = c \in \mathbb{N}$, we obtain the class of sequences $n_F=c\cdot n$ for $n\geq 1$. Naturally layers of such cobweb posets designated by the sequence satisfying (\ref{eq:1}) for $n\geq 1$  may also be partitioned according to (cprta1).
\vspace{0.4cm}

\noindent \textbf{Example 1.5} $1_F = 1, n_F = c\cdot n, n>1$ ($0_F$ = corresponds to one "\emph{empty root}" $\{\emptyset\}$ )
 This might be considered another sample example now illustrating the "\emph{shifted}" method named (cpta2). For example if we choose $2_F = c \in \mathbb{N}$, while $1_F=1$, we obtain the class of sequences $1_F = 1$ and $n_F = c\cdot n $ for $n>1$. Layers of such cobweb posets designated by these sequences may also be partitioned.

\begin{observen} \emph{\textbf{ Algorithm (cpta2) }}
\noindent Given any (including cobweb-admissible) sequence $A\equiv\{n_A\}_{n\geq 0}$, $s\in\mathbb{N}\cup\{0\}$   let us define \emph{shift} unary operation $\oplus_s$ as follows:
\begin{displaymath}
	\oplus_s A = B, \qquad n_B = \left\{ 
	\begin{array}{ll} 
		1 				& n<s \\ 
		(n-s)_A 		& n\geq s   
	\end{array} \right. 
\end{displaymath}

\noindent where $B\equiv\{n_B\}_{n\geq 0}$. Naturally $\oplus_0$ = identity.  Then the following is true.
If a sequence $A$ is cobweb-tiling sequence then $B$ is also cobweb-tiling sequence.
\end{observen}

\vspace{0.4cm}
\noindent For example this is the case for $A=1,2,3,4,\ldots$, $\oplus_{3}A = 1,1,1,1,2,3,4,\ldots$.

\begin{example}
$n_F = m_F \cdot k_F$
\end{example} 
If we choose $1_F = c\in \mathbb{N}$, we obtain the class of sequences $n_F = c^n, n\geq0$. We can also consider more general case $n_F = \alpha\cdot m_F \cdot k_F$, where $\alpha\in\mathbb{N}$ which gives us the next class of tiling sequences $n_F = \alpha^{n-1}\cdot c^n, n\geq1, 0_F=1$ and layers of such cobweb posets can be partitioned by (cprta1) algorithm. For example: $1_F = 1, \alpha=2 \rightarrow F = 1,1,2,4,8,16,32,\ldots$ or $1_F = 2 \rightarrow F = 1,2,4\alpha,8\alpha^2,16\alpha^3,\ldots$

\begin{example}
$n_F = (m+k)_F = (k+1)_F\cdot m_F + (m-1)_F\cdot k_F$
\end{example}
Here also we have infinite number of cobweb tiling sequences depending on the initial values chosen for the recurrence $(k\!+\!2)_F\!=\!2_F(k\!+\!1)_F\!+\!k_F, k\!\!\geq\!\!0$. For example: $1_F = 1$ and $2_F = 2 \rightarrow$ $F=1,2,5,12,29,70,169,408,985,\ldots$ 
Note that this is not shifted Fibonacci sequence as we use recurrence (2)  which depends on initial conditions adopted.
Next $1_F = 1$ and $2_F = 3 \rightarrow F=1,3,10,33,109,360,1189,\ldots$ 
Note that this is not remarcable Lucas sequence [7].

\vspace{0.4cm}
Neither of sequences: shifted Fibonacci nor Lucas sequence  satisfy (2)
neither these (as well as the Catalan, Motzkin, Bell or Euler numbers sequences)  are cobweb admissible sequences.  This indicates the further exceptionality   of Fibonacci sequence   along with natural numbers.

\vspace{0.4cm}
The proof of tiling existence leads to many easy known formulas for sequences, where we use multiplications of terms $m_F$ and/or $k_F$, like
$n_F = \alpha\cdot k_F$,
$n_F = \alpha\cdot m_F k_F$,
$n_F = \alpha\cdot (m\pm\beta)_F k_F$,
where $\alpha, \beta \in \mathbb{N}$, $n=m+k$ and so on.

This are due to the fact that in the course of partition's existence proving with (cprta1) partition of layer $\langle\Phi_{k+1}\!\rightarrow\!\Phi_n\rangle$ existence relies on partition's existence of smaller layers $\langle\Phi_{k+1}\!\rightarrow\!\Phi_{n-1}\rangle$ and/or $\langle\Phi_{k}\!\rightarrow\!\Phi_{n-1}\rangle$.

\vspace{0.4cm}

In what follows we shall use an at the point product of two cobweb-admissible sequences giving as a result a new cobweb admissible sequence - cobweb tiling sequences included to which the above described treatment (cprta1) applies.

\subsection{Beginnings of the cobweb-admissible sequences$\\ $ production}

\begin{defn}
Given any two cobweb-admissible sequences $A\equiv\{n_A\}_{n\geq 0}$ and $B\equiv\{n_B\}_{n\geq 0}$, their at the point product $C$ is given by
\begin{displaymath}
	A\cdot B = C \qquad C\equiv\{n_C\}_{n\geq 0},\ n_C = n_A \cdot n_B
\end{displaymath}
\end{defn}

\noindent It is obvious  that $A \cdot B = C$ is also cobweb admissible and 
\begin{displaymath}
	{n \choose k}_{A\cdot B} = \frac{n_{A}^{\underline{k}}}{k_A!} \cdot \frac{n_{B}^{\underline{k}}}{k_B!} = {n \choose k}_A \cdot {n \choose k}_B \in \mathbb{N}\cup\{0\}
\end{displaymath}

\begin{example}
Almost constant sequences $C_t$ \label{ex:const}
\begin{displaymath}
	C_t = \{n_C\}_{n\geq 0}\qquad \mathrm{where}\  const = n_C = t \in \mathbb{N}\ \mathrm{for}\ n>0, 0_F = 1.
\end{displaymath}
\end{example}

\noindent as for example $C_5 = 1,5,5,5,5,\ldots$ are trivially cobweb-admissible and cobweb tiling sequences - see next example. 

\vspace{0.4cm}
\noindent In the following $I$ denotes unit sequence $I\equiv\{1\}_{n\geq 0}$; $I\cdot{A}=A$.

\begin{example}
Not diminishing sequence $A_{c,M} $
\end{example}

\noindent If we multiply $i$-th term (where $i\geq M \geq 1, M \in \mathbb{N}$) of sequence $I$ by any constant $c \in \mathbb{N}$, then the product cobweb admissible sequence is $A_{c,M}$.
\begin{displaymath}
	A_{c,M} \equiv \{n_A\}_{n\geq 0}\qquad \mathrm{where}\ n_A = 
	\left\{ 
	\begin{array}{lr} 
		1 		& 1\leq n < M \\ 
		c 		& n\geq M   
	\end{array} \right. 
\end{displaymath}

\noindent as for example $\nolinebreak{A_{5,10} = 1,\underbrace{1,\ldots,1}_{10},5,5,5,\ldots}$ or more general example  

\noindent $\nolinebreak{A_{3,2,10} = 1,\underbrace{3,\ldots,3}_{10},6,6,6,\ldots}$ Clearly sequences of this type are cobweb admissible and  cobweb tiling sequences.

\vspace{0.4cm}
Indeed. Each of level of layer $\langle\Phi_k\!\rightarrow\!\Phi_n\rangle$ has the same or more vertices than each of levels of the block $\sigma P_m$. If not the same then the number of vertices from the block $\sigma P_m$  divides the number of vertices at corresponding layer's level. This is how   (cprta2) applies.

\vspace{0.4cm}
\noindent \textbf{Note.}
The  sequence $A_{3,2,10}$ is a product of two sequences from Example \ref{ex:const}, 
$A = 1,3,3,3,3,3,3,\ldots$ and  $B'=\oplus_{10}B= 1,\ldots,1,2,2,2,\ldots$ where $\\ B = 1,2,2,2,2,2,2,\ldots$, then $A\cdot B' = A_{3,2,10} = 1,\underbrace{3,\ldots,3}_{10},6,6,6,\ldots $ 

\begin{example}
Periodic sequence $B_{c,M}$
\end{example}

\noindent A more general example is supplied by 
\begin{displaymath}
	B_{c,M} \equiv \{n_B\}_{n\geq 0}\qquad \mathrm{where}\ n_B = 
	\left\{ 
	\begin{array}{ll} 
		1 		& M\nmid{n} \vee n=0 \\ 
		c 		& M | n  
	\end{array} \right. 
\end{displaymath}

\noindent where $c,M \in \mathbb{N}$. Sequences of above form are cobweb tiling, as for example  $B_{2,3} = \underbrace{1,1,2}_{3},1,1,2,\ldots$, $B_{7,4} = \underbrace{1,1,1,7}_{4},1,1,1,7,\ldots$ Indeed.

% % PROOF
\vspace{0.4cm}
\noindent \textbf{P}\textbf{\footnotesize{ROOF}}
\vspace{0.2cm}
\noindent Consider any layer $\langle\Phi_k\!\rightarrow\!\Phi_n\rangle$, $k\leq n$, $k,n\in\mathbb{N}\cup\{0\}$, with $m$ levels:

\noindent For $m < M$, the block $P_m$ has one vertex on each of levels. The tiling is trivial.

\noindent For $m \geq K$, the sequence $B_{c,M}$ has a period equal to $M$, therefore any layer of $m$ levels has the same or larger number of levels with $c$ vertices than the block $\sigma P_m$, if layer's level has more vertices than corresponding level of block $\sigma P_m$ then the quotient of this numbers is a natural number i.e. $1|c$, thus the layer can be partitioned by one block $P_m$ or by $c$ blocks $\sigma P_m$ $\blacksquare$

\begin{observen}
$\\ $\emph{
The at the point product of the above sequences gives us occasionally a method to produce Natural numbers as well as expectedly other cobweb-admissible sequences with help of the following algorithm.}
\end{observen}

\vspace{0.4cm}
\noindent \textbf{Algorithm for natural numbers' generation (cta3)}

\vspace{0.4cm}
\noindent $N(s)$ denotes a sequence which first $s$ members is next Natural numbers i.e. 
$N(s)\equiv\{n_N\}_{n\geq 0}$, where $n_N = n$, for $n=1,2,\ldots,s$, $\ p, p_n$ - prime numbers.

\begin{enumerate}
	\item $N(1) = \mathbf{I} = 1,1,1,\ldots $
	\item $N(2) = N(1)\cdot B_{2,2} = 1,2,1,2,1,2,\ldots$
	\item $N(3) = N(2)\cdot B_{3,3} = 1,2,3,2,1,6,\ldots$
	\item[n.] $N(n) = N(n-1) \cdot \mathbf{X} $
\end{enumerate}

\noindent Consider $n$:
\begin{enumerate}
	\item let $n$ be prime, then $\neg \exists_{1\neq i \in [n-1]} i|n \Rightarrow n_N = 1 \Rightarrow \mathbf{X} = B_{n,n} $

	\item let $n = p^m$, $1<m \in \mathbb{N}$, then $n_N = p^{m-1} \Rightarrow \mathbf{X} = B_{p,n}$

	\item let $n = \prod_{s=1}^{u} p_{s}^{m_s}$, where $p_i \neq p_j$ for $i\neq j$, $m_i \geq 1$, $\\ i=1,2,\ldots,u$, $u>1$
	
	$\forall_{i\in[u]} p_i^{m_i} < n  \Rightarrow n_N = \mathrm{LCD}\left(\{ p_i^{m_i}:  i=1,2,\ldots,u \}\right) $
	$\\ \wedge \forall_{i\neq j} \mathrm{GCD}( p_i^{m_i},  p_j^{m_j}) = 1 $ 
	$\Rightarrow$ $n_N = \prod_{s=1}^{u} p_{s}^{m_s}$
	$\Rightarrow \mathbf{X} = \mathbf{I}$

\end{enumerate}

\noindent where lowest common denominator or least common denominator (LCD)  and  greatest common divisor (GCD)  abbreviations were used.

\vspace{0.4cm}
\noindent \textbf{Concluding}
\begin{displaymath}
	N(n) = N(n-1)\cdot B_{h_n, n} 
	\begin{array}{c}
		\scriptstyle{n \rightarrow \infty} \\
		\longrightarrow \\
		\ 
	\end{array}
	\mathbb{N}
\end{displaymath}
\begin{displaymath}	
	h_n = \left\{ 
	\begin{array}{ll}
	p & \quad n = p^m,\quad \mathbb{N} \ni m \geq 1  \\
	1 & \quad n = \prod_{s=1}^{u>1} p_{s}^{m_s},\quad \mathbb{N} \ni m_s \geq 1
	\end{array}
	\right.
\end{displaymath}

\vspace{0.2cm} \noindent 
while  $\{h_n\}_{n\geq 1} = 1,2,3,2,5,1,7,2,3,1,11,1,13,1,1,2,17,\ldots$

\begin{figure}[ht]
\begin{center}
	\includegraphics[width=100mm]{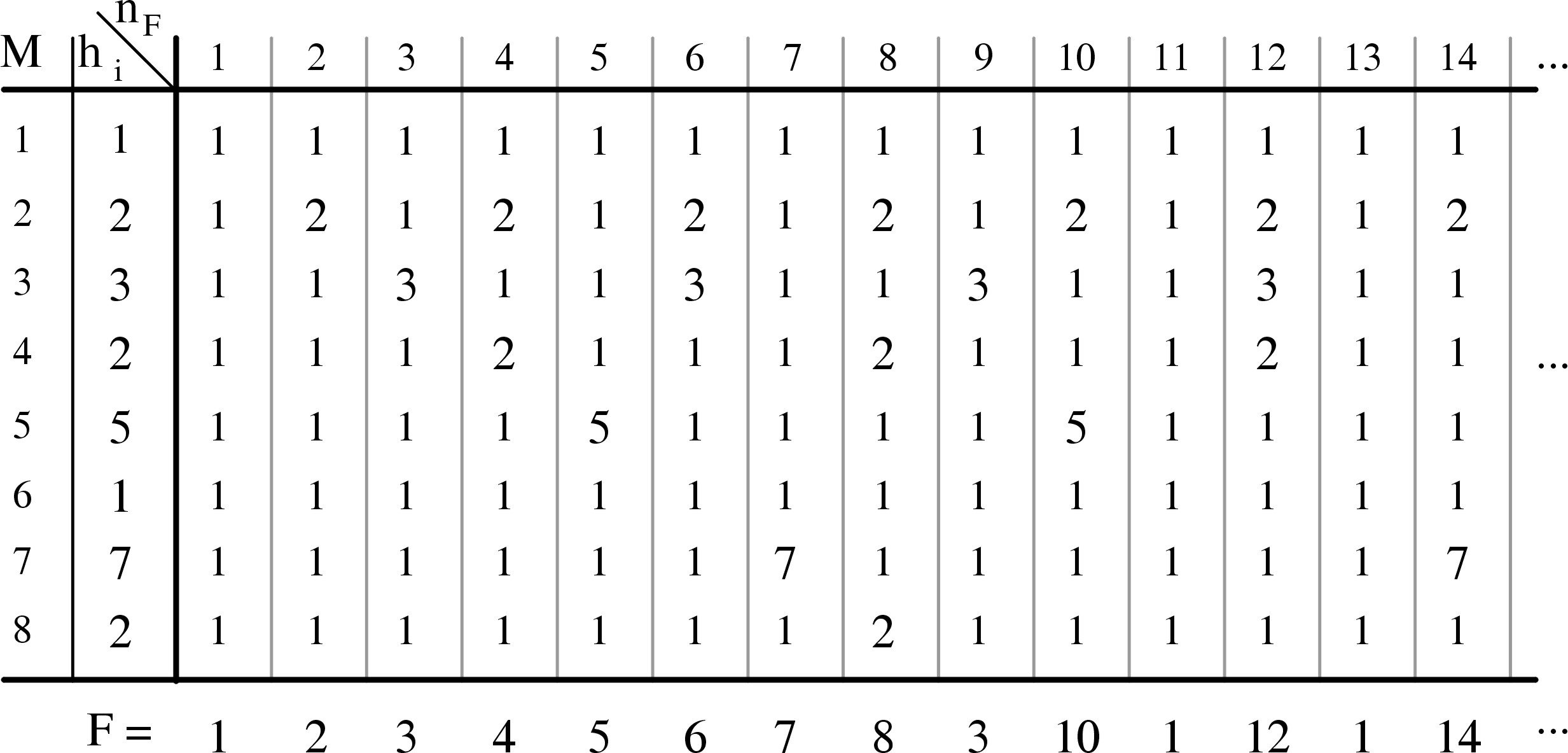}
	\caption{Display of eight steeps of algorithm (cta3) \label{fig:gennat}}
\end{center}
\end{figure}

\vspace{0.2cm}
As for the Fibonacci sequence we expect the same statement to be true for $n\rightarrow\infty$ bearing in mind those properties of  Fibonacci numbers which make them an effective tool in Zeckendorf representation of natural numbers. For the Fibonacci numbers the would be sequence $\{h_n\}_{n\geq 1}$  is given by $\{h_n\}_{n\geq 1}=1,1,2,3,5,4,13,7,17,11,89,6,\ldots$

\vspace{0.2cm}
We end up with general observation - rather obvious but important to be noted.

\begin{theoremn}
Not all cobweb-admissible sequences are cobweb tiling sequences. 
\end{theoremn}

\noindent \textbf{P}\textbf{\footnotesize{ROOF}}

\begin{figure}[ht]
\begin{center}
	\includegraphics[width=110mm]{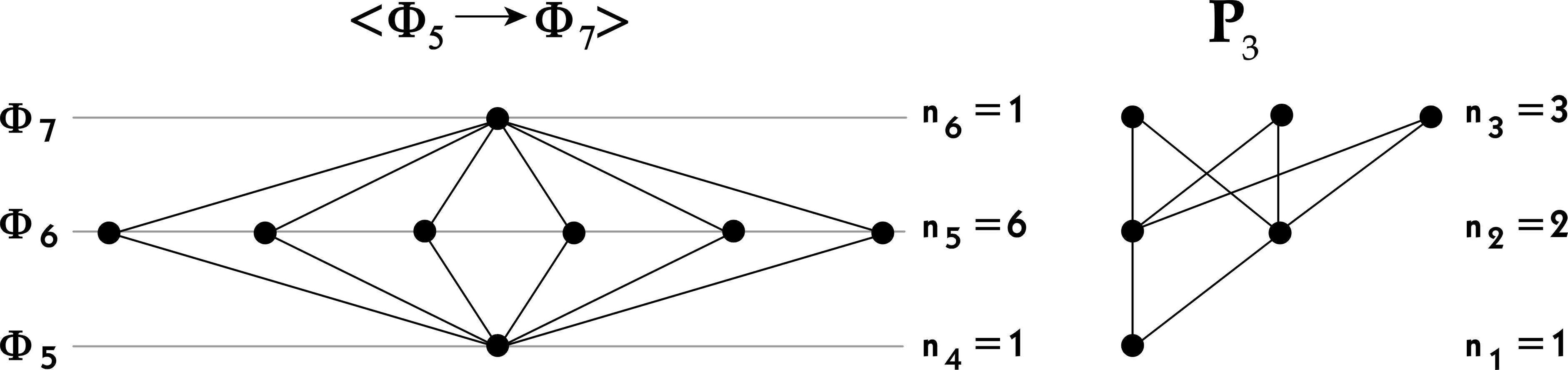}
	\caption{Picture proof of Theorem 3 \label{fig:contr} }
\end{center}
\end{figure}

\noindent It is enough to give an appropriate example. Consider then a cobweb-admissible sequence $F=A\cdot B=1,2,3,2,1,6,1,2,3,\ldots$, where $A=1,2,1,2,1,2\ldots$  and $B=1,1,3,1,1,3,\ldots$ are both cobweb admissible and cobweb tiling.
Then the layer $\langle\Phi_5\!\rightarrow\!\Phi_7\rangle$ can not be partitioned with blocks $\sigma P_3$ as the level $\Phi_5$ has one vertex, level $\Phi_5$ has six while $\Phi_5$ has one vertex again (Fig \ref{fig:contr}).

\vspace{0.4cm}
\noindent \textbf{Corollary}
\emph{The at the point product of two tiling sequences does not need to be a tiling sequence.}

\vspace{0.4cm}
However  for  $A=1,2,1,2,\ldots$ and $B=1,1,3,1,1,3,\ldots$ cobweb tiling sequences their product $F=A\cdot B=1,2,3,2,1,6,1,\ldots$ is not a cobweb  tiling sequence. A natural question - is it still ahead [6,7]? . Find the effective characterizations and or algorithms for a cobweb admissible sequence to be a cobweb tiling sequence.  The second author has encoded the problem of an algorithm being looked for with help of his invention
called by him "`Primary cobweb admissible binary tree"'  and this is a subject of a separate note to be presented soon.

\section{GCD-morphism Problem. Problem III.} 

\noindent \textbf{Coming over to the last problem announced} above following  [6-8]  let us note that the Observation 4. provides us
with the \textit{new} combinatorial interpretation of  the immense  class
of all classical $F-nomial$ coefficients including
binomial or Gauss $q$- binomial ones or Konvalina generalized binomial
coefficients of the first and of the second kind [3] - which include Stirling numbers of both kinds too.
All these  $F$-nomial coefficients naturally are computed with their correspondent
cobweb-admissible sequences. More than that - the vast `umbral' family of  
$F$-sequences [9-13,4] includes also those which are called \textit{ "GCD-morphic"} sequences.
This means that $GCD[F_n,F_m] = F_{GCD[n,m]}$ where $GCD$ stays for Greatest Common Divisor. 
\begin{defn}.
The sequence of integers $F=\{n_F\}_{n\geq 0}$  is called the
GCD-morphic  sequence  if  $GCD[F_n,F_m] = F_{GCD[n,m]}$ where
$GCD$ stays for Greatest Common Divisor operator.
\end{defn}
The Fibonacci sequence is a much nontrivial [11,12,6] guiding example of GCD-morphic sequence.
Of course  \textit{not all }incidence coefficients of reduced incidence
algebra of full binomial type are computed with  GCD-morphic sequences
however these or that - if computed with the cobweb correspondent admissible sequences  
all are given the new, joint  cobweb poset combinatorial interpretation via Observation 3. 
More than that - in [8] a prefab-like combinatorial description of cobweb posets is being served with corresponding
generalization of the fundamental exponential formula. 

\noindent \textbf{Question:} which of these above mentioned sequences are GCD-morphic
sequences?

\vspace{2mm}

\noindent \textbf{GCD-morphism Problem. Problem III.} \textit{ Find effective characterizations
and/or  an algorithm to produce the GCD-morphic sequences i.e. find all examples.}

\noindent The second author has "`almost solved"' the GCD-morphism Problem - again  with help of his invention
called by him "`Primary cobweb admissible binary tree"'  and this is a subject of a separate note to be presented soon.(See [27]).

\vspace{2mm}

\section {Apendix}

\noindent \textbf{A.1. Cobweb posets and  KoDAGs' ponderables  of Kwa\'sniewski relevant recent productions.} [19-25,7,6]

\begin{defn}
\noindent  Let  $n\in N \cup \left\{0\right\}\cup \left\{\infty\right\}$. Let   $r,s \in N \cup \left\{0\right\}$.  Let  $\Pi_n$ be the graded partial ordered set (poset) i.e. $\Pi_n = (\Phi_n,\leq)= ( \bigcup_{k=0}^n \Phi_k ,\leq)$ and $\left\langle \Phi_k \right\rangle_{k=0}^n$ constitutes ordered partition of $\Pi_n$. A graded poset   $\Pi_n$  with finite set of minimal 
elements is called \textbf{cobweb poset} \textsl{iff}  
$$\forall x,y \in \Phi \  i.e. \  x \in \Phi_r \ and \  y \in \Phi_s \   r \neq s\ \Rightarrow \   x\leq y   \ or \ y\leq x  , $$ 
 $\Pi_\infty \equiv \Pi. $
\end{defn}

\vspace{0.1cm}

\noindent \textbf{Note}. By definition of $\Pi$ being graded its  levels    $\Phi_r \in \left\{\Phi_k\right\}_k^\infty$ are independence sets  and of course partial order  $\leq $ up there in Definition 6.1. might be replaced by $<$.

\vspace{0.2cm}

\noindent The Definition 11 is the reason for calling Hasse digraph $D = \left\langle \Phi, \leq \cdot \right\rangle $ of the poset $(\Phi,\leq))$ a \textbf{\textcolor{red}{Ko}}DAG as in  Professor   
\textbf{\textcolor{red}{K}}azimierz   \textbf{\textcolor{red}{K}}uratowski native language one word \textbf{\textcolor{red}{Ko}mplet} means \textbf{complete ensemble} - see more in [19-25].
 
\begin{defn}
\noindent Let  $F = \left\langle k_F \right\rangle_{k=0}^n$ be an arbitrary natural numbers valued sequence, where $n\in N \cup \left\{0\right\}\cup \left\{\infty\right\}$. We say that the cobweb poset $\Pi = (\Phi,\leq)$ is \textcolor{red}{\textbf{denominated}} (encoded=labelled) by  $F$  iff   $\left|\Phi_k\right| = k_F$ for $k = 0,1,..., n.$
\end{defn}

\vspace{0.2cm}

\noindent \textbf{A.2.}  See also much relevant  [26,2011]

\vspace{0.4cm}

\noindent \textbf{A.3. Cobweb posets and  combinatorial interpretation in discrete hyper-boxes language.}[19],[29]

\vspace{0.1cm}

\noindent \textbf{Theorem}. [19] \\
\noindent \textit{For $F$-cobweb admissible sequences $F$-binomial coefficient} $\fnomial{n}{k}$ \textit{is the cardinality of the family of \emph{equipotent} to}  $V_{0,m}$ \textit{mutually disjoint 
discrete hyper-boxes, all together \textbf{partitioning } the discrete hyper-box }  $V_{k+1,n}$    $\equiv$   \textit{the layer}   $\langle\Phi_{k+1} \to \Phi_n \rangle$,\textit{ where} $m=n-k$.

\vspace{0.3cm}

\noindent  \textbf{\textcolor{blue}{The cobweb tiling problem} in the language of \textcolor{blue}{discrete hyper-boxes}}.

\vspace{0.2cm}

%%%%%%%%%%%%%%%%%%%%%%%%%%%%%%%%%%%%%%%%%%%%%%%%%%%%%%%%%%%%%%%%%%%%%%%%%%%%%%%%%
\noindent \textbf{Comment} General  "fractal-reminiscent" comment. The discrete $m$-dimensional $F$-box ($m = n-k$) with edges' sizes designated by natural numbers' valued sequence $F$ 
where invented in [26] as a response to the so called \emph{ cobweb tiling problem} posed in [6,2007] and then repeated in [7,2009]. This tiling problem was considered by Maciej Dziemia\'nczuk in [1,2008] where it was shown that not all admissible $F$-sequences permit tiling as defined in [6,2007].  Then - after [26,2009 ArXiv] this tiling  problem was considered by  Maciej Dziemia\'nczuk in discrete hyper-boxes language \cite[2009]{48}.

\noindent \textbf{Recall the fact} ([6,2007], [7,2009]): \textit{Let} $F$\textit{ be an admissible sequence}. \textit{Take any natural numbers} $n,m$ such that $n\geq m$, \textit{then the value of} $F$-\textit{binomial coefficient} $\fnomial{n}{k}$  \textit{s equal to the number} of sub-boxes \textit{that constitute a}  $\kappa$-\textit{partition of} $m$-\textit{\textit{dimensional}} $F$-\textit{box} $V_{m,n}$ \textit{where} $\kappa = |V_m|$.

\begin{defn}
Let $V_{m,n}$ be a $m$-dimensional $F$-box. Then any $\kappa$-partition into sub-boxes of the form $V_m$ is called tiling of $V_{m,n}$.
\end{defn}

\noindent Hence \textcolor{blue}{\textbf{only these}} partitions of $m$-dimensional box $V_{m,n}$  are admitted  for which all sub-boxes \textbf{are of the form} $V_m$ i.e. we have a kind of  (\textcolor{blue}{\textbf{self-similarity}}).\\

\noindent It was shown in \cite[2008]{13} by Maciej Dziemia\'nczuk  that the only admissibility condition  is not sufficient for the existence a tiling for any given $m$-dimensional box $V_{k,n}$. Kwa\'sniewski in [6,2007] and [7,2009] posed the question called \emph{Cobweb Tiling Problem} which we repeat here.

\vspace{0.1cm}

\noindent \textbf{Tiling problem}\\
\noindent Suppose that $F$ is an admissible sequence. Under which conditions any $F$-box $V_{m,n}$ designated by sequence $F$ has a tiling? Find effective characterizations and/or find an algorithm to produce these tilings.

\vspace{0.1cm}

\noindent In \cite[2009]{48} by Maciej Dziemia\'nczuk one  proves the existence of such tiling for certain sub-family of admissible sequences $F$. These  include among others $F=$ Natural numbers, Fibonacci numbers, or $F =  \left\langle n_q \right\rangle_{n\geq 0}$ Gaussian sequence. Original extension of the above  tiling problem onto the general case multi $F$-multinomial coefficients is  proposed  in \cite[2009]{48} , too. Moreover - a reformulation of the present cobweb tiling problem into a clique problem of a graph specially invented for that purpose - is invented.

%%%%%%%%%%%%%%%%%%%%%%%%%%%%%%%5

\vspace{2mm}
\begin
{thebibliography}{99}
\parskip 0pt

\bibitem{1}
Maciej  Dziemia\'nczuk, {\it On Cobweb posets tiling problem}, Adv. Stud. Contemp. Math. volume \textbf{16} (2), (\textbf{2008}): 219-233. arXiv:0709.4263v2, Thu, 4 Oct 2007 14:13:44 GMT

\bibitem{2}
Charles Jordan {\it On Stirling Numbers}  Tôhoku Math. J. \textbf{37} (\textbf{1933}),254-278.     

\bibitem{3}
John Konvalina , {\it A Unified Interpretation of the Binomial Coefficients, the Stirling Numbers and the Gaussian Coefficients}
The American Mathematical Monthly {\bf 107} (2000), 901-910.

\bibitem{4}
Ewa Krot, {\it An Introduction to Finite Fibonomial Calculus}, CEJM 2(5) (2005) 754-766.

\bibitem{5}
Ewa Krot, {\it The first ascent into the Fibonacci Cob-web Poset},  Adv. Stud. Contemp. Math.  \textbf{11} (2) (\textbf{2007}) 179-184.

\bibitem{6}
A. Krzysztof Kwa\'sniewski, \emph{Cobweb posets as noncommutative prefabs}, Adv. Stud. Contemp. Math.   vol.14 (1) (2007): 37 - 47. arXiv:math/0503286v4 [v4] Sun, 25 Sep 2005 23:40:37 GMT

\bibitem{7}
A. Krzysztof  Kwa\'sniewski, \textit{On cobweb posets and their combinatorially admissible sequences}, Advanced Studies in Contemporary Mathematics, \textbf{18} no. 1,  (\textbf{2009}): 17-32.  arXiv:math/0512578v5,  [v5] Mon, 19 Jan 2009 21:47:32 GMT

\bibitem{8}
A. Krzysztof Kwa\'sniewski, {\it First observations on Prefab posets' Whitney numbers},  Advances in Applied Clifford Algebras Volume 18, Number 1 / February, 2008, 57-73, Xiv:0802.1696v1, [v1] Tue, 12 Feb 2008 19:47:18 GMT

\bibitem{9}
A. Krzysztof  Kwa\'sniewski,{\it On extended finite operator calculus of Rota and quantum groups} Integral Transforms and Special Functions, {\bf 2} (4) (2001) 333-340.

\bibitem{10}
A. Krzysztof  Kwa\'sniewski, {\it Main  theorems of extended finite operator calculus} Integral Transforms and Special Functions, {\bf 14} (6) (\textbf{2003}) 499-516.

\bibitem{11}
A. Krzysztof  Kwa\'sniewski, {\it The  Logarithmic Fib-Binomial Formula}, Adv. Stud. Contemp. Math. v.9 No.1 (\textbf{2004}): 19 -26. arXiv:math/0406258v1 [v1] Sun, 13 Jun 2004 17:24:54 GMT

\bibitem{12}
A. Krzysztof  Kwa\'sniewski, {\it Fibonomial cumulative connection constants}, Bulletin of the ICA  \textbf{44} (\textbf{2005}) 81- 92. Upgrade: arXiv:math/0406006v6, [v6] Fri, 20 Feb 2009 02:26:21

\bibitem{13}
A. Krzysztof  Kwa\'sniewski, {\it On umbral extensions of Stirling numbers and Dobinski-like formulas}, Advanced Stud. Contemp. Math. {\bf 12}(\textbf{2006}) no. 1, pp.73-100.
arXiv:math/0411002v5, [v5] Thu, 20 Oct \textbf{2005} 02:12:47 GMT

\bibitem{14}
Anatoly D. Plotnikov, {\it About presentation of a digraph by dim 2 poset}, Adv. Stud. Contemp. Math.  \textbf{12} (1) (2006) 55-60

\bibitem{15}

Bruce E. Sagan {\it Mobius Functions of Posets} (Lisbon lectures)IV: \textit{Why the Characteristic Polynomial factors }June 28 \textbf{2007} 
http://www.math.msu.edu/%7Esagan/Slides/mfp4.pdf

\bibitem{16}
E. Spiegel, Ch. J. O`Donnell  {\it Incidence algebras}  Marcel
Dekker, Inc., Basel, 1997.

\bibitem{17}

Richard P. Stanley, \textit{Hyperplane Arrangements}, Proc. Nat.  Acad. Sci. \textbf{93} (1996), 2620-2625.
\textit{An Introduction to Hyperplane Arrangements } www.math.umn.edu/~ezra/PCMI2004/stanley.pdf

\bibitem{18}
Morgan Ward: {\em A calculus of sequences}, Amer.J.Math. \textbf{58} (1936)  255-266.

\bibitem{19}                                                                                                                                                              
A. Krzysztof Kwa\'sniewski \textit{Natural join construction of graded posets versus ordinal sum and discrete hyper boxes}, arXiv:0907.2595v2 [v2] Thu, 30 Jul 2009 22:46:39 GMT

\bibitem{20}                                                                                                                                                                
A. Krzysztof Kwasniewski \textit{On natural join of posets properties and first applications} arXiv:0908.1375v2 [v2]  Sat, 22 Aug 2009 10:42:44 GMT   

\bibitem{21}                                                                                                                                                          
A. Krzysztof Kwa\'sniewski \textit{Graded posets inverse zeta matrix  formula}, Bull. Soc. Sci. Lett. Lodz , vol. 60,  no.3  (2010)  in print.  arXiv:0903.2575v3 [v3] Mon, 24 Aug 2009 05:38:45 GMT

\bibitem{22}                                                                                                                                                              
A. Krzysztof Kwa\'sniewski \textit{Graded posets zeta matrix formula}, Bull. Soc. Sci. Lett. Lodz , vol. 60,  no. 3 (2010),   in print .  arXiv:0901.0155v2 Thu, Mon, 16 Mar 2009 15:43:08 GMT

\bibitem{23}                                                                                                                                                         
A.K. Kwa\'sniewski , {\it Some Cobweb Posets Digraphs' Elementary Properties and Questions}, Bull. Soc. Sci. Lett. Lodz , vol. 60,  no. 2  (2010) in print. arXiv:0812.4319v1 [v1] Tue, 23 Dec 2008 00:40:41 GMT

\bibitem{24}                                                                                                                                                             
A. Krzysztof Kwa\'sniewski , {\it Cobweb Posets and  KoDAG  Digraphs are Representing Natural Join of Relations, their di-Bigraphs and the Corresponding  Adjacency Matrices}, Bull. Soc. Sci. Lett. Lodz , vol. 60,  no. 1 (2010)  in print.  arXiv:0812.4066v3,[v3] Sat, 15 Aug 2009 05:12:22 GMT 

\bibitem{25}                                                                                                                                                          
A. Krzysztof Kwa\'sniewski, \emph{How the work of Gian Carlo Rota had influenced my group research and life}, arXiv:0901.2571v4 Tue, 10  Feb 2009 03:42:43 GMT

\bibitem{26} 
A. Krzysztof  Kwa\'sniewski, Maciej  Dziemia\'nczuk , {\it On cobweb posets' most relevant codings}, to appear in Bull. Soc. Sci. Lett. Lodz , vol. 61 (2011), arXiv:0804.1728v2 [v2] Fri, 27 Feb 2009 18:05:33 GMT

\bibitem{27} 
Maciej Dziemia\'nczuk, Wies{\l}aw, Wladys{\l}aw Bajguz, \emph{On GCD-morphic sequences}, IeJNART: Volume (\textbf{3}), September (\textbf{2009}): 33-37. arXiv:0802.1303v1, [v1] Sun, 10 Feb 2008 05:03:40 GMT 

\bibitem{48}                                                                                                                                                                              
Maciej  Dziemia´nczuk, \textit{On cobweb posets and discrete F-boxes tilings}, arXiv:0802.3473v2  v2,  Thu, 2 Apr\textbf{ 2009}.11:05:55 GMT      

\bibitem{29}
A. Krzysztof  Kwa\'sniewski, \textit{Note on Ward-Horadam H(x) - binomials' recurrences and related interpretations, II}, to appear; (scheduled to be announced at Mon, 10 Jan 2011). 

\end{thebibliography}
 \end{document}